\documentclass[11pt,fullpage]{article}

\usepackage[pdftex]{graphicx}
\usepackage{latexsym}
\usepackage{color}
\usepackage{multirow}
\usepackage{amsmath}
\usepackage{amssymb}
\usepackage{amsthm}
\usepackage{amsfonts}
\usepackage{bm}

\usepackage[normalem]{ulem}

\usepackage{xr,refcount}

\makeatletter
\newcommand*{\addFileDependency}[1]{
  \typeout{(#1)}
  \@addtofilelist{#1}
  \IfFileExists{#1}{}{\typeout{No file #1.}}
}
\makeatother

\addFileDependency{SupplementaryMaterial.tex}
\addFileDependency{SupplementaryMaterial.aux}

\newtheorem{theorem}{Theorem}
\newtheorem{lemma}{Lemma}

\newtheorem*{assumption}{Assumption}
\newtheorem{corollary}{Corollary}

\long\def\ignore#1{}

\newcommand{\be}{\begin{equation}}
\newcommand{\ee}{\end{equation}}
\newcommand{\beqn}{\begin{eqnarray}}
\newcommand{\eeqn}{\end{eqnarray}}

\newcommand{\bfm}[1]{\mbox{\boldmath{$#1$}}}
\newcommand{\bbeta}{\bfm{\beta}}
\newcommand{\btheta}{\bfm{\theta}}

\newcommand{\bxi}{\bfm{\xi}}

\newcommand{\bp}{{\bf p}}

\newcommand{\bx}{{\bf x}}
\newcommand{\bX}{{\bf X}}

\newcommand{\cB}{{\cal B}}

\newcommand{\cb}{{\cal B}}

\newcommand{\cE}{{\cal E}}
\newcommand{\cF}{{\cal F}}
\newcommand{\cC}{{\cal C}}

\newcommand{\cX}{{\cal X}}
\newcommand{\cH}{{\cal H}}

\newcommand{\M}{\mathfrak{M}}

\newcommand{\sparsity}{\Delta}
\newcommand{\complexity}{r}
\renewcommand{\top}{T}

\begin{document}

\title{\bf Multiclass classification by sparse multinomial logistic regression}

\author{{\bf Felix Abramovich}\\
Department of Statistics\\
 and Operations Research \\
Tel Aviv University \\
Israel \\
{\it felix@tauex.tau.ac.il}
\and
{\bf Vadim Grinshtein} \\
Department of Mathematics\\ 
and Computer Science\\
The Open University of Israel\\ 
Israel\\
{\it vadimg@openu.ac.il}
\and
{\bf Tomer Levy} \\
Department of Statistics\\
 and Operations Research \\
Tel Aviv University \\
Israel \\
{\it tmrlvi@gmail.com}
}

\date{}

\maketitle

\begin{abstract}
In this paper we consider high-dimensional multiclass classification by sparse multinomial logistic regression. 
We propose first a feature selection procedure based on penalized maximum likelihood with
a complexity penalty on the model size and derive the nonasymptotic bounds for misclassification excess risk of the resulting classifier. We establish also their tightness by deriving the corresponding minimax lower bounds. In particular, we show that there exist two regimes corresponding to small and large number of classes.
The bounds can be reduced under the additional low noise condition. 
To find a  penalized maximum likelihood solution with a complexity penalty requires, however, a combinatorial search over
all possible models. To design a feature selection procedure computationally feasible for high-dimensional data, we propose 
multinomial logistic group Lasso and Slope classifiers and show that they also achieve the minimax order.
\end{abstract}


\bigskip

\section{Introduction} \label{sec:intr}
Classification is one of the core problems in statistical learning and has been intensively studied in statistical and machine learning literature.  
Nevertheless, while the theory for binary classification is 
well developed (see, Devroy, Gy\"ofri and Lugosi, 1996; Vapnik, 2000; Boucheron, Bousquet and Lugosi, 2005 and references therein for a comprehensive review), its multiclass extensions are much less complete. 

Consider a general $L$-class classification with a (high-dimensional) vector of features $\bX \in \cX \subseteq \mathbb{R}^d$ and the outcome class label $Y \in \{1,\ldots,L\}$. We can model it as
$Y|(\bX=\bx) \sim Mult(p_1(\bx),\ldots,p_L(\bx))$, where
$p_l(\bx)=P(Y=l|\bX=\bx),\;l=1,\ldots,L$.

A classifier is a measurable function $\eta: \cX \rightarrow \{1,\ldots,L\}$. The accuracy of a classifier $\eta$
is defined by a misclassification error $R(\eta)=P(Y \neq \eta(\bx))$. 
The optimal classifier that minimizes this error is
the Bayes classifier $\eta^*(\bx)=\arg \max_{1 \leq l \leq L}p_l(\bx)$ with $R(\eta^*)=1-E_{\bX}\max_{1 \leq l \leq L}p_l(\bx)$.
The probabilities $p_l(\bx)$'s are, however, unknown and one should derive a classifier $\widehat{\eta}(\bx)$ from the available data $D$: a random sample of $n$ independent observations $(\bX_1,Y_1),\ldots, (\bX_n,Y_n)$ from the joint distribution of $(\bX,Y)$. 
The corresponding (conditional) misclassification error of $\widehat{\eta}$ is $R(\widehat{\eta})=P(Y \neq \widehat{\eta}(\bx)|D)$ and the goodness of $\widehat{\eta}$ w.r.t. $\eta^*$ is measured by the misclassification
excess risk $\cE(\widehat{\eta},\eta^*)=ER(\widehat{\eta})-R(\eta^*)$. The goal is then to find a classifier $\widehat{\eta}$ within given family with minimal $\cE(\widehat{\eta},\eta^*)$.

A first strategy in multiclass classification is to reduce it to a series of binary classifications. The probably two most well-known methods
are One-vs-All (OvA), where each class is compared against all others, and One-vs-One (OvO), where all 
pairs of classes are compared to each other.  

A more direct and appealing strategy is to extend binary classification approaches for multiclass case. 
Thus, a common approach to design a multiclass classifier $\widehat{\eta}$ is based on empirical risk minimization (ERM), where minimization of 
a true misclassification error $R(\eta)$  is replaced by minimization of the corresponding empirical risk
$\widehat{R}_n(\eta)=\frac{1}{n}\sum_{i=1}^n I\{Y_i \neq \eta(\bx_i)\}$ over a given class of classifiers. For binary classification, tight risk bounds
for ERM classifiers have been established in terms of 
VC-dimension, Rademacher complexity or covering numbers 
(see Devroy, Gy\"ofri and Lugosi, 1996; Vapnik, 2000; Boucheron, Bousquet and Lugosi, 2005 and references therein). Their extensions to multiclass case,
however, are not straightforward. See  Maximov and Reshetova (2016) for a comprehensive survey of the state-of-the-art results on the upper bounds for misclassification excess risk of multiclass ERM
classifiers.
A comparison of error bounds for ERM classifiers with those for
OvA and OvO is given in Daniely {\em et al.} (2012).

A crucial drawback of ERM is in minimization of 0-1 loss that makes it computationally infeasible. A typical remedy is
to replace 0-1 loss by a related convex surrogate.  The resulting solution approximates then the minimizer of the corresponding surrogate risk. The goal is to find a surrogate loss such that minimization of its risk leads to a Bayes classifier $\eta^*$ (aka Fisher consistent or calibrated loss). Various calibrated losses for multiclass classification have been considered in
the literature (e.g., Zhang, 2004b; Chen and Sun, 2006; Tewari and Bartlett, 2007; 
\'Avila Pires, Szepesv\'ari and Ghavamzadeh, 2013;
\'Avila Pires and Szepesv\'ari, 2016).

An alternative approach to ERM is to estimate $p_l(\bx)$'s from the data by some $\widehat{p}_l(\bx)$'s
and to use a plug-in classifier of the form  $\widehat{\eta}(\bx)=
\arg \max_{1 \leq l \leq L} \widehat{p}_l(\bx)$.
A standard approach is to assume some (parametric or nonparametric) model for $p_l(\bx)$. 
The most commonly used model is multinomial logistic regression,
where it is assumed
that
$p_l(\bx)=\frac{\exp(\bbeta^T_l \bx)}{\sum_{k=1}^L \exp(\bbeta^T_k \bx)}$ and $\bbeta_l \in \mathbb{R}^d,\;l=1,\ldots,L$ are unknown vectors of regression coefficients. The corresponding Bayes classifier 
is, therefore, a linear classifier $\eta^*(\bx)=\arg \max_{1 \leq l \leq L}p_l(\bx)=\arg \max_{1 \leq l \leq L} \bbeta^T_l \bx$.
One then estimates $\bbeta$'s from the data by the maximum likelihood estimators (MLE) $\widehat{\bbeta}_l$'s and derives
the plug-in (linear) classifier $\widehat{\eta}(\bx)=\arg \max_{1 \leq l \leq L} \widehat{\bbeta}^T_l \bx$.  Unlike ERM, the MLE  $\widehat{\bbeta}$'s though not available in the closed form, can be nevertheless obtained numerically by the fast iteratively reweighted least squares
algorithm (McCullagh and Nelder, 1989, Section 2.5).

The general challenge modern statistics faces with is high-dimensionality of the data, where the number of features $d$ is large and might be
even larger than the sample size $n$ (large $d$ small $n$ setups) that raises a severe
``curse of dimensionality'' problem. 
Reducing the dimensionality of a feature space by selecting a sparse subset of ``significant'' features becomes crucial.

For binary classification Devroy, Gy\"orfi and Lugosi (1996, Chapter 18) and Vapnik (2000, Chapter 4) considered model selection from a sequence of classifiers within
a sequence of classes  by penalized ERM with the structural penalty depending
on the VC-dimension of a class. See also  Boucheron, Bousquet and Lugosi (2005, Section 8) 
for related penalized ERM approaches and references therein. 
Abramovich and Grinshtein (2019) explored feature selection in high-dimensional logistic regression classification. 

To the best of our knowledge, feature selection for multiclass classification has not yet been rigorously well-studied and the goal of this paper is to fill the gap.
Thus, we propose a model/feature selection procedure 
based on penalized maximum likelihood
with a certain complexity penalty on the model size. We establish
the non-asymptotic upper bounds for misclassification excess risk of the resulting plug-in
classifier which is also adaptive to the unknown sparsity and show their tightness  by deriving the corresponding
minimax lower bound over a set of sparse linear classifiers. It turns out that there appear interesting phenomena in the multiclass setup. 
In particular, we find that there exist two regimes. For $L \leq 2+\ln(d/d_0)$, where $d_0$ is
the size of the true (unknown) model, the multiclass
effect is not manifested and
the minimax misclassification excess risk over the set of $d_0$-sparse linear classifiers is of the order $\sqrt{\frac{d_0}{n}\ln\left(\frac{de}{d_0}\right)}$ regardless of $L$. For larger $L$, it increases as $\sqrt{\frac{d_0 (L-1)}{n}}$ and does not depend on $d$. 
We also
show that these bounds can be improved under the additional low-noise assumption.

Any penalized maximum likelihood procedure that involves a complexity
penalty requires, however, a combinatorial search over all
possible models that makes its use computationally infeasible for large $d$.
A common remedy is then to use a convex surrogate, where the original combinatorial minimization is replaced by a related convex program. In this paper
we consider Slope convex relaxation which can be viewed as
generalization of the celebrated Lasso and show that for the properly chosen tuning parameters, the resulting multinomial logistic group Slope multiclass classifier is also minimax rate-optimal.

The rest of the paper is organized as follows. In Section \ref{sec:class} we present sparse multinomial logistic regression model and propose a feature selection procedure. The 
bounds for misclassification
excess risk of the resulting plug-in classifier are derived in
Section \ref{sec:bounds}. 
In Section \ref{sec:low-noise} we introduce the additional low-noise assumption that allows one to improve the bounds. In Section \ref{sec:convex} we develop group
Slope convex relaxation techniques for multiclass classification with Lasso as its particular case, and establish the misclassification excess risk bounds for the resulting classifier. All the proofs are given in the Appendix.

\section{Construction of a classifier} \label{sec:class}
\subsection{Multinomial logistic regression model} \label{subsec:multinom}
Consider $d$-dimensional $L$-class classification model that can be written in the following form: 
\be \label{eq:model}
Y|(\bX=\bx) \sim Mult(p_1(\bx),\ldots,p_L(\bx)), 
\ee
where $\bX \in \mathbb{R}^d$ is a vector of linearly independent features with a marginal probability distribution $P_X$ with a support $\cX \subseteq \mathbb{R}^d$ and $\sum_{j=1}^L
p_j(\bx)=1$ for any $\bx \in \cX$. 

We consider a multinomial logistic regression model, where it is assumed that
\be \label{eq:multinom}
\ln \frac{p_l(\bx)}{p_L(\bx)}=\bbeta^T_l \bx,\;\;\;l=1,\ldots,L-1,
\ee
and $\bbeta_l \in \mathbb{R}^d$ are the vectors of the (unknown)
regression coefficients.
Hence,
$$
p_l(\bx)=\frac{\exp(\bbeta_l^T \bx)}{1+\sum_{k=1}^{L-1}\exp(\bbeta_k^T \bx)},\;\;l=1,\ldots,L-1\;\;\;{\rm and}\;\;\;
p_L(\bx)=\frac{1}{1+\sum_{k=1}^{L-1}\exp(\bbeta_k^T \bx)},
$$
or, in a somewhat more compact form,
$$
p_l(\bx)=\frac{\exp(\bbeta_l^T \bx)}{\sum_{k=1}^L\exp(\bbeta_k^T \bx)},\;\;l=1,\ldots,L
$$
with $\bbeta_L={\bf 0}$.
We set $\bbeta_l={\bf \infty}=(\infty,\ldots,\infty)$ and $\bbeta_l={\bf -\infty}=(-\infty,\ldots,-\infty)$ to include two degenerate cases $p_l(\bx)=1$ and $p_l(\bx)=0$ respectively.

The Bayes classifier is then a linear classifier 
$\eta^*(\bx)=\arg \max_{1 \leq l \leq L} p_l(\bx)=
\arg \max_{1 \leq l \leq L} \bbeta^T_l \bx$ with misclassification
risk $R(\eta^*)=1-E_\bX\max_{1 \leq l \leq L} p_l(\bx)$.

The choice of the last class as a reference class is, in fact, quite arbitrary.
One can consider an equivalent model with any other reference class $h$ instead:
$\ln \frac{p_l(\bx)}{p_h(\bx)}=\gamma^T_l \bx,\;\;\;l \neq h$. Evidently, there is one-to-one transformation: $\gamma_l=\bbeta_l-\bbeta_h$ and
$\bbeta_l=\gamma_l-\gamma_L$. Change of a reference class is, therefore, just a matter of reparametrization of the same model. 

\subsection{Penalized maximum likelihood estimation} \label{subsec:mle}
To each possible value $y \in \{1,\ldots,L\}$ of $Y$ assign the  indicator
vector $\bxi \in \{0,1\}^L$ with $\xi_l=I\{y=l\},\;l=1,\ldots,L$.
Let $B \in \mathbb{R}^{d \times L}$ be the matrix of the regression
coefficients in (\ref{eq:multinom}) with the columns $\bbeta_1,\ldots,\bbeta_L$ (recall that $\bbeta_L={\bf 0}$) and 
let $f_B(\bx,y)$ be the corresponding joint 
distribution of $(\bX,Y)$, i.e.
$d f_B(\bx,y)=\prod_{l=1}^L p_l(\bx)^{\xi_l}~ dP_X(\bx)$,
where $p_l(\bx)=\frac{\exp\{\bbeta^T_l \bx \}}{\sum_{k=1}^L \exp\{\bbeta^T_k \bx \}}$.
Given a random sample $(\bX_1,Y_1),\ldots,(\bX_n,Y_n) \sim f_B(\bX,Y)$, the conditional log-likelihood function is
\be \label{eq:likelihood}
\ell(B)=\sum_{i=1}^n \left\{\bX_i^T B \bxi_i-\ln \sum_{l=1}^L\exp(\bbeta_l^T \bX_i)\right\},
\ee
and one can find the maximum likelihood estimator (MLE) for $B$ by maximizing $\ell(B)$.

The era of ``Big Data'' brought the challenge of dealing with problems, where the number of features $d$ 
is very large and may be even larger than the sample size $n$ (``large $d$ small $n$''
setups). Nevertheless, it is commonly assumed that the true underlying model is sparse and most of the features do not have a significant impact on classification. Reducing the dimensionality of a feature space by selecting a sparse subset of ``significant'' features is then crucial. Thus, Bickel and Levina (2004) and Fan and Fan (2008) showed that even binary high-dimensional classification without a proper feature selection might be as bad as just pure guessing.

For binary classification, where the regression matrix $B$ reduces to a single vector $\bbeta \in \mathbb{R}^d$, the sparsity is naturally measured by the $l_0$ (quasi)-norm $|\bbeta|_0$ -- the number of non-zero entries of $\bbeta$ (see, e.g., Abramovich and Grinshtein, 2019).  
For multiclass case one can think of several possible ways to extend the 
notion of sparsity. The most evident measure of sparsity is the number of non-zero rows of $B$ that corresponds to the assumption that part of the features do not have any impact on classification at all and, therefore, have zero coefficients in (\ref{eq:multinom}) for all $l$. It can be viewed as {\em global} row-wise sparsity.  
One can easily verify that such a measure is invariant under the choice of the
reference class in (\ref{eq:multinom}). 

In what follows we assume the following assumption:
\begin{assumption}[\bf A] \label{as:A}
Assume that there exists $0 < \delta < 1/2$ such that $\delta < 
p_l(\bx) < 1-\delta$
or, equivalently, $|\bbeta_l^T \bx| < C_0$ with $C_0=\ln \frac{1-\delta}{\delta}$ for all $\bx \in \cX$ and all $l=1,\ldots,L$. 
\ignore{
In particular, for $|\bx|_2 \leq 1$, it is sufficient to assume the
boundedness of $\bbeta_l$'s, i.e.
$|B|_{2,\infty}=\max_{1 \leq l \leq L}|\bbeta_l|_2 \leq C_0$.
}
\end{assumption}
Assumption (A) prevents the conditional variances $Var(\xi_l|\bX=\bx)=p_l(\bx)(1-p_l(\bx))$ to be infinitely close to zero, where any MLE-based
procedure may fail. 

Let $\M$ be the set of all $2^d$ possible models $M \subseteq
\{1,\ldots,d\}$. 
In view of Assumption (A), for a given model $M$ define a set of matrices $\cb_M =\{B \in \mathbb{R}^{d \times L}: B_{\cdot L}={\bf 0},\;|B|_{2,\infty} \leq \ln\frac{1-\delta}{\delta} \;{\rm and}\;B_{j \cdot}={\bf 0}\;\;{\rm iff}\; j \not\in M\}$. Obviously, all matrices in $\cb_M$ have the same number of non-zero rows which can be naturally defined as a model
size $|M|$.

Under the model $M$, the MLE $\widehat{B}_M$ of $B$ 
is then
\be \label{eq:mle}
\widehat{B}_M=\arg \max_{\widetilde{B} \in \cb_M}
\sum_{i=1}^n \left\{\bX_i^T \widetilde{B} \bxi_i-\ln \sum_{l=1}^L\exp(\widetilde{\bbeta}_l^T \bX_i)\right\},
\ee
where $\widetilde{\bbeta}_l=\widetilde{B}_{\cdot l},\;l=1,\ldots L$ are the columns of $\widetilde{B}$.
\ignore{
The corresponding (linear) plug-in classifier 
$$
\eta_M(\bx)=\arg \max_{1 \leq l \leq L} \widehat{p}_{Ml}(\bx)=
\arg \max_{1 \leq l \leq L} \widehat{\bbeta}_{Ml}^T \bx
$$
}

Select the model $\widehat{M}$ by the penalized maximum likelihood model selection criterion of the form 
\be \label{eq:binpen}
\widehat{M}=
\arg \min_{M \in \M}\left\{ 
\sum_{i=1}^n \left(\ln\left(\sum_{l=1}^L\exp(\widehat{\bbeta}_{Ml}^T \bX_i)\right)-\bX_i^T \widehat{B}_M \bxi_i\right)+Pen(|M|) \right\}
\ee 
with the complexity penalty $Pen(\cdot)$ on the model size $|M|$.

Finally, for the selected model $\widehat{M}$ the resulting plug-in classifier
\be \label{eq:class}
\widehat{\eta}_{\widehat M}(\bx)=\arg \max_{1 \leq l  \leq L} \widehat{\bbeta}_{\widehat M l}^T \bx
\ee

The proper choice of the complexity penalty $Pen(\cdot)$ in
(\ref{eq:binpen}) is obviously the core of the proposed approach.

\section{Misclassification excess risk bounds} \label{sec:bounds}
We now derive the (non-asymptotic) upper bound for misclassification excess risk of
the penalized maximum likelihood classifier (\ref{eq:class})
derived in Section \ref{sec:class} for a particular type of the
complexity penalty and then show that such a choice is, in fact, optimal (in the minimax sense).

Denote the number of nonzero rows of a matrix $B$ by $r_B$.
Let $\cC_L(d_0)=\{\eta(\bx)=\arg \max_{1 \leq l \leq L} \bbeta^T_l\bx: B \in \mathbb{R}^{d \times L},\; B_{\cdot L}={\bf 0}\;{\rm and}\;
r_B \leq d_0\}$ be 
the set of all $d_0$-sparse linear $L$-class classifiers. The sparsity parameter $d_0$ is assumed to be unknown and the goal is to construct classifiers adaptive to the unknown sparsity. 

\begin{theorem} \label{th:upper}
Consider a $d_0$-sparse multinomial logistic regression
model (\ref{eq:model})-(\ref{eq:multinom}). 

Let $\widehat{M}$ be a model selected in
(\ref{eq:mle})-(\ref{eq:binpen}) with the complexity penalty
\be \label{eq:optpen} 
Pen(|M|)=c_1 |M|(L-1)+c_2 |M| \ln\left(\frac{d e}{|M|}\right),
\ee 
where the absolute constants $c_1, c_2 >0$ are given in the proof of Theorem \ref{th:upper_low-noise}. 

Then, under Assumption (A), 
\be \label{eq:upper}
\sup_{\eta^* \in \cC_L(d_0)}\cE(\widehat{\eta}_{\widehat M},\eta^*) \leq  C_1(\delta)~\sqrt{\frac{d_0 (L-1)+   
d_0 \ln \left(\frac{de}{d_0}\right)}{n}}
\ee
for some $C_1(\delta)$ depending on $\delta$, simultaneously for all $1 \leq d_0 \leq \min(d,n)$.
\end{theorem}
Theorem \ref{th:upper} is a particular case of a more general
Theorem \ref{th:upper_low-noise} given in the next Section \ref{sec:low-noise}.

The complexity penalty $Pen(|M|)$ in (\ref{eq:optpen}) contains two terms.
The first one is proportional to $|M|(L-1)$ -- the overall number
of estimated parameters in the model $M$ and is an AIC-type
penalty. The second one is proportional to $|M|
\ln \left(\frac{de}{|M|}\right) \sim  \ln \binom{d}{|M|}$ -- the log(number of all possible models of size $|M|$) and typically appears in model selection in various regression and classification setups (see, e.g.
Birg\'e and Massart, 2001; Bunea, Tsybakov and Wegkamp, 2007; Abramovich and Grinshtein, 2010, 2016, 2019). 

Theorem \ref{th:lower} below shows that for an agnostic model,
where the Bayes risk $R(\eta^*)>0$, the upper bound (\ref{eq:upper}) for the misclassification excess risk established in Theorem \ref{th:upper} is essentially tight and up to a possibly different constant coincides
with the corresponding minimax lower bound over $\cC_L(d_0)$:
\begin{theorem} \label{th:lower}
Consider a $d_0$-sparse agnostic multinomial logistic regression
model (\ref{eq:model})-(\ref{eq:multinom}), where
$2 \leq  d_0 \ln(\frac{de}{d_0}) \leq n$ and $d_0(L-1) \leq n$. Then,
$$
\inf_{\widetilde \eta}\sup_{\eta^* \in \cC_L(d_0),~P_X}
\cE(\widetilde{\eta},\eta^*) \geq C_2 {\sqrt \frac{d_0 (L-1)+d_0
\ln \left(\frac{de}{d_0}\right)}{n}}
$$
for some $C_2>0$, where the infimum is taken over all classifiers $\tilde{\eta}$ based on the data $(\bX_i,Y_i),\;i=1,\ldots,n$.
\end{theorem}

The above bounds imply, in particular, that as $d$ and $L$ grow with $n$ and assuming that $\delta$ is bounded away from zero, there are two different regimes: 

\medskip
\noindent
{\em 1. Small number of classes: $L \leq 2+\ln\left(\frac{d}{d_0}\right)$.}
\newline
In this case, the complexity penalty (\ref{eq:optpen}) is
$Pen(|M|) \sim c |M| \ln\left(\frac{d e}{|M|}\right)$ does not depend on $L$. The resulting 
(rate-optimal) misclassification excess risk is of the order
$\sqrt{\frac{d_0}{n} \ln \left(\frac{de}{d_0}\right)}$ regardless of $L$ and the error in feature selection dominates 
in the overall excess risk.
Multiclass classification for such a small number of classes  is essentially not harder than binary
(see the results of Abramovich and Grinshtein, 2019 for $L=2$). 

\medskip
\noindent
{\em 2. Large number of classes: $2+\ln\left(\frac{d}{d_0}\right) < L \leq \frac{n}{d_0} $.}
\newline
In this regime, $Pen(|M|) \sim c |M|(L-1)$ is an AIC type penalty (see above), the misclassification excess risk increases with $L$ as $\sqrt{\frac{d_0(L-1)}{n}}$ regardless of $d$ and the main contribution to the overall error comes from estimating the large number ($d_0(L-1)$) of parameters in the model.

\medskip
For $L > \frac{n}{d_0}$ the number of parameters in the model becomes
larger than the sample size and consistent classification is
impossible.

In particular, without sparsity assumption, i.e. in the case $d_0=d~ (\leq n)$, the misclassification excess risk is of the order $\sqrt{\frac{d(L-1)}{n}}$ for
all $1 \leq L-1 \leq \frac{n}{d}$.

Note that even if the considered multinomial logistic regression model is misspecified and the Bayes
classifier $\eta^*$ is not linear, we still have the
following risk decomposition
\be \label{eq:decomp}
R(\widehat{\eta}_{\widehat M})-R(\eta^*)=\left(R(\widehat{\eta}_{\widehat M})-R(\eta^*_L)\right)+
\left(R(\eta^*_L)-R(\eta^*)\right),
\ee
where $\eta^*_L=\arg \min_{\eta \in \cC_L(d)} R(\eta)$ is the oracle (ideal) linear classifier.  
The above results can then be still applied to the first term in the RHS of (\ref{eq:decomp})
representing the estimation error, while the second term is an approximation
error and measures the ability of linear classifiers to perform as good as $\eta^*$.
Enriching the class of classifiers may improve
the approximation error but will necessarily increase the estimation error in (\ref{eq:decomp}). 
In a way, it is similar to the variance/bias tradeoff in regression.

\section{Improved bounds under low-noise condition} \label{sec:low-noise}
Intuitively, it is clear that misclassification error is particularly large when it is difficult to separate the class with the highest probability from others, i.e. at those $\bx \in \cX$, where the 
largest probability $p_{(1)}(\bx)$ is close to the second largest $p_{(2)}(\bx)$ (see also
Kesten and Morse, 1959). 

Define the following multiclass extension of the low-noise (aka Tsybakov) condition (Mammen and Tsybakov, 1999; Tsybakov, 2004):
\begin{assumption}[\bf B] \label{as:B}
Consider the multinomial logistic regression model (\ref{eq:model})-(\ref{eq:multinom}) and assume that there exist
$C>0, \alpha \geq 0$ and $h^*>0$ such that for all $0<h \leq h^*$,
\be \label{eq:low-noise}
P\left(p_{(1)}(\bX)-p_{(2)}(\bX) \leq h\right) \leq C h^{\alpha} 
\ee
\end{assumption}
\noindent
(see also Chen and Sun, 2006; Agarwal, 2013).
Assumption (B) implies that with high probability (depending on the parameter $\alpha$) the most likely  class is sufficiently distinguished from others. The two extreme cases $\alpha=0$ and
$\alpha=\infty$ correspond respectively to no assumption on the noise considered in the previous sections and to existence of a hard margin of size $h^*$ separating $p_{(1)}(\bx)$ and $p_{(2)}(\bx)$. 
A straightforward multiclass extension of Lemma 5 of Bartlett, Jordan and
McAuliffe (2006) implies that (\ref{eq:low-noise}) is equivalent to
the condition that there exists $C_1(\alpha)$ such that for any
classifier $\eta$,
\be \label{eq:low-noise_equiv}
P\left(\eta(\bX) \neq \eta^*(\bX)\right) \leq  C_1(\alpha)~ 
\cE(\eta,\eta^*)^{\frac{\alpha}{\alpha+1}}
\ee

We now show that under the additional low-noise condition (\ref{eq:low-noise}) the bounds for the
misclassification excess risks established in the previous Section \ref{sec:bounds} can be improved:
\begin{theorem} \label{th:upper_low-noise}
Consider a $d_0$-sparse multinomial logistic regression model (\ref{eq:model})-(\ref{eq:multinom}) and
let $\widehat{M}$ be a model selected in
(\ref{eq:binpen}) with the complexity penalty (\ref{eq:optpen}).

Then, under Assumptions (A) and (B), there exists $C(\delta)$ such that
$$
\sup_{\eta^* \in \cC_L(d_0)}\cE(\widehat{\eta}_{\widehat M},\eta^*) \leq C(\delta) \left(\frac{d_0 (L-1)+   
d_0 \ln \left(\frac{de}{d_0}\right)}{n}\right)^{\frac{\alpha+1}{\alpha+2}}
$$
for all $1 \leq d_0 \leq \min(d,n)$ and all $\alpha \geq 0$.
\end{theorem}
Thus, $\widehat{\eta}_{\widehat M}$ is adaptive to both
$d_0$ and $\alpha$.
As we have mentioned, Theorem \ref{th:upper} is a particular case of Theorem \ref{th:upper_low-noise} with $\alpha=0$. 

To conclude this section we note that the error bounds can be also improved under other types of additional constraints on the marginal distribution $P_X$, e.g, a so-called strong density assumption (Audibert and Tsybakov, 2007 for binary classification) or a cluster assumption  (Rigollet, 2007; Maximov, Amini and Harchaoui, 2018).

\section{Multinomial logistic group Lasso and Slope} \label{sec:convex}
Despite strong theoretical results on penalized maximum likelihood classifiers with complexity penalties  established in 
the previous sections, selecting the model $\widehat{M}$ in (\ref{eq:binpen}) requires a combinatorial search over all possible models in $\M$ that makes it computationally infeasible when the
number of features is large. A common approach to handle this problem is convex relaxation, where the original combinatorial minimization is replaced by a related convex surrogate. The most well-known examples 
include the celebrated Lasso, where the $l_0$-norm in the penalty is replaced by $l_1$-norm
norm, and its recently developed more general variation Slope that uses a {\em sorted} $l_1$-type norm (Bogdan {\em et al.}, 2015). Lasso and Slope estimators have been intensively studied in
Gaussian regression (see, e.g., 
Bickel, Ritov and Tsybakov, 2009; Su and Cand\'es, 2015;
Bellec, Lecu\'e and Tsybakov, 2018 among others), and their logistic modifications in
logistic regression (van de Geer, 2008; Abramovich and Grinshtein, 2019; Alquier, Cottet and Lecu\'e, 2019). Abramovich and Grinshtein (2019) investigated logistic Lasso and Slope classifiers for the binary case. In this section we consider {\em multinomial} logistic {\em group} Lasso and Slope classifiers and extend the corresponding results of Abramovich and Grinshtein (2019) for multiclass classification.

Recall that we consider a 
global row-wise sparsity, where the coefficient regression matrix $B$ has a subset of zero rows. 
To capture such type of sparsity we consider a {\em multinomial} logistic
{\em group} Lasso and Slope classifiers defined as follows. For a given tuning parameter $\lambda>0$, find
\be \label{eq:groupLasso}
\widehat{B}_{gL} =\arg \min_{\widetilde B}
\left\{\frac{1}{n}\sum_{i=1}^n\left(\ln\left(\sum_{l=1}^L \exp(\widetilde{\bbeta}_l^T \bX_i)\right)-\bX_i^T\widetilde{B}\xi_i\right)+\lambda \sum_{j=1}^d |{\widetilde B}|_j \right\},
\ee
where $|\widetilde{B}|_j=|\widetilde{B}_{j\cdot}|_2$ is the $l_2$-norm of the $j$-th row of $\widetilde{B}$ and define the corresponding classifier 
$\widehat{\eta}_{gL}(\bx)=\arg\max_{1 \leq l \leq L}\widehat{\bbeta}_{gL,l}^T \bx$.
An efficient algorithm for computing
multinomial logistic group Lasso is given in Vincent and Hansen (2014).

Multinomial logistic group Slope is a more general variation of (\ref{eq:groupLasso}). Namely, 
\be \label{eq:groupSlope}
\widehat{B}_{gS} =\arg \min_{\widetilde B}
\left\{\frac{1}{n}\sum_{i=1}^n\left(\ln\left(\sum_{l=1}^L \exp(\widetilde{\bbeta}_l^T \bX_i)\right)-\bX_i^T\widetilde{B}\bxi_i\right)+ \sum_{j=1}^d \lambda_j |{\widetilde B}|_{(j)} \right\},
\ee
where the rows' $l_2$-norms  $|\widetilde{B}|_{(1)} \geq \ldots \geq
|\widetilde{B}|_{(d)}$ are the descendingly ordered
and $\lambda_1 \geq \ldots \geq \lambda_d > 0$  are the tuning parameters, and set $\widehat{\eta}_{gS}(\bx)=\arg\max_{1 \leq l \leq L}\widehat{\bbeta}_{gS,l}^T \bx$. 
Multinomial logistic group Lasso (\ref{eq:groupLasso}) can be evidently viewed as a particular case of (\ref{eq:groupSlope}) with equal $\lambda_j$'s. 

Note that unlike complexity penalties, the solution of (\ref{eq:groupSlope}) is identifiable without the additional constraint $\widetilde{\bbeta}_L={\bf 0}$.
Moreover, since the unconstrained log-likelihood (\ref{eq:likelihood}) satisfies $\ell(\widetilde{\bbeta}_1,\ldots,\widetilde{\bbeta}_L)=\ell(\widetilde{\bbeta}_1-{\bf c},\ldots,\widetilde{\bbeta}_L-{\bf c})$ for any  
vector ${\bf c} \in \mathbb{R}^d$, the solution can be always improved by taking $\hat{c}_j=\arg\min_{c_j} \sum_{l=1}^L (\widetilde{B}_{jl}-c_j)^2$, i.e. 
$\hat{c}_j=\bar{B}_{j\cdot}$. Hence, $\widehat{B}_{gS}$ necessarily satisfies the symmetric constraint
$\sum_{l=1}^L \widehat{\bbeta}_{gS,l}={\bf 0}$ or, equivalently, $\widehat{B}_{gS} {\bf 1}={\bf 0}$ (zero mean rows).

As usual for convex relaxation methods, one needs some (mild) constraints
on the design. In particular, we assume the
following assumption on the marginal distribution $P_X$:
\begin{assumption}[\bf C] Assume that for all (generally dependent) components $X_j$'s of a random features vector $\bX \in \mathbb{R}^d$,
\begin{enumerate}
\item $EX_j^2=1$ ($X_j$'s are scaled)
\item there exist constants $\kappa_1, \kappa_2, w > 1$ and $\gamma \geq 1/2$ such that $E(|X_j|^p)^{1/p} \leq \kappa_1 p^{\gamma}$ for all $2 \leq p \leq \kappa_2 \ln(w d)$
($X_j$'s have polynomially growing moments up to the order $\ln d$)
\end{enumerate}
\end{assumption}
In particular, Assumption (C) evidently holds for (scaled) Gaussian and sub-Gaussian $X_j$'s with $\gamma=1/2$ for all moments.  Assumption (C) ensures that
for $n \geq C_1(\ln d)^{\max(2\gamma-1,1)}$, 
\be \label{eq:lm17}
E\max_{1 \leq j \leq d} \frac{1}{n} \sum_{i=1}^n X_{ij}^2 \leq C_2
\ee
for some constants $C_1=C_1(\kappa_1,\kappa_2,w,\gamma)$ and $C_2=C_2(\kappa_1,\kappa_2,w)$ (Lecu\'e and Mendelson, 2017,
proof of Theorem A). Moreover, (\ref{eq:lm17}) might be violated if the moments condition in Assumption (C) holds only up to
the order of $\ln(wd)/\ln\ln(wd)$.  We will need (\ref{eq:lm17}) in the proof of the upper bound for misclassification excess risk of a general multinomial logistic group Slope classifier (\ref{eq:groupSlope}).

For simplicity of exposition, in what follows we consider
$\gamma=1/2$ corresponding to $n \geq C_1 \ln d$, where $C_1$ is given in
Lecu\'e and Mendelson (2017).

\begin{theorem} \label{th:slope} Consider a $d_0$-sparse multinomial logistic regression (\ref{eq:model})-(\ref{eq:multinom}). 

Apply the multinomial logistic group Slope classifier (\ref{eq:groupSlope}) with 
$\lambda_j$'s satisfying
\be \label{eq:lambda_slope}
\max_{1 \leq  j \leq d} \frac{\sqrt{L+\ln(d/j)}}{\lambda_j} \leq C_0\sqrt{n}
\ee
with the constant $C_0$ derived in the proof.

Assume Assumptions (A)-(C) and let $n \geq C_1 \ln d$.

Then, 
$$
\sup_{\eta^* \in \cC_L(d_0)}\cE(\widehat{\eta}_{gS},\eta^*) 
\le C(\delta) \left(\sum_{j=1}^{d_0} 
\frac{\lambda_j}{\sqrt{j}}\right)^{\frac{2(\alpha+1)}{\alpha+2}}   
$$
for some constant $C(\delta)$ depending on $\delta$.
\end{theorem}

We now consider two specific choices of $\lambda_j$'s: 

\medskip
\noindent
1. {\em Equal $\lambda_j$ (multinomial logistic group Lasso)}.

Take
\be \label{eq:lambda_lasso}
\lambda = C_0 \sqrt{\frac{L +\ln d}{n}} 
\ee
to satisfy (\ref{eq:lambda_slope}).
Note that
$\sum_{j=1}^{d_0} \frac{1}{\sqrt j} \leq 2 \sqrt{d_0}$ that yields the following corollary of Theorem \ref{th:slope}:
\begin{corollary} \label{cor:lasso}
Consider a $d_0$-sparse multinomial logistic regression (\ref{eq:model})-(\ref{eq:multinom}).
Apply the multinomial logistic group Lasso classifier (\ref{eq:groupLasso}) with $\lambda$ from (\ref{eq:lambda_lasso}).

Then, under Assumptions (A)-(C) and  $n \geq C_1 \ln d$,
$$
\sup_{\eta^* \in \cC_L(d_0)}\cE(\widehat{\eta}_{gL},\eta^*) \leq C(\delta) \left(\frac{d_0 (L-1) +   
d_0 \ln(de)}{n}\right)^{\frac{\alpha+1}{\alpha+2}}
$$
for all $1 \leq d_0 \leq \min(d,n)$ and all $\alpha \geq 0$.
\end{corollary}
Thus, unless $d$ grows faster than exponentially with $n$,
the multinomial logistic group Lasso classifier $\widehat{\eta}_{gL}$ achieves a minimax order for large number of classes (see Section \ref{sec:bounds}), while for small $L$ it is rate-optimal for sparse cases, where $d_0 \ll d$, but only sub-optimal (up to an extra logarithmic loss) for dense cases, where $d_0 \sim d$.  
We conjecture that similar to the results of Bellec, Lecu\'e and Tsybakov (2018) for Gaussian regression, $\widehat{\eta}_{gL}$ with {\em adaptively} chosen $\lambda$ can achieve the minimax rate in the latter case as well but the proof of this conjecture is beyond the scope of the paper.

\medskip
\noindent
{\em 2. Variable $\lambda_j$'s.} Consider  
\be \label{eq:lambda_slope0}
\lambda_j = C_0 \sqrt{\frac{L+\ln(d/j)}{n}}
\ee
that evidently satisfies (\ref{eq:lambda_slope}).
One can also verify that 
$$
\sum_{j=1}^{d_0} \sqrt{\frac{L+\ln(d/j)}{j}} \leq \frac{2L}{L-1} \sqrt{d_0(L+\ln(d/d_0))} \leq 4\sqrt{d_0\left(L-1+\ln\left(\frac{de}{d_0}\right)\right)}
$$  
Theorem \ref{th:slope} implies then:

\begin{corollary} \label{cor:slope} 
Consider a $d_0$-sparse multinomial logistic regression (\ref{eq:model})-(\ref{eq:multinom}).
Apply the multinomial logistic group Slope classifier (\ref{eq:groupSlope}) with $\lambda_j$'s from (\ref{eq:lambda_slope0}).

Then, under Assumptions (A)-(C) and  $n \geq C_1 \ln d$,
$$
\sup_{\eta^* \in \cC_L(d_0)}\cE(\widehat{\eta}_{gS},\eta^*) \leq C(\delta) \left(\frac{d_0 (L-1) +   
d_0 \ln\left(\frac{de}{d_0}\right)}{n}\right)^{\frac{\alpha+1}{\alpha+2}}
$$
for all $1 \leq d_0 \leq \min(d,n)$ and all $\alpha \geq 0$.
\end{corollary}
Hence, if the number of features grows at most exponentially with $n$,
the multinomial logistic group Slope classifier with variable $\lambda_j$'s from (\ref{eq:lambda_slope0}) is adaptively rate-optimal for both small and large number of classes, and, unlike the penalized likelihood classifier $\widehat{\eta}_{\widehat M}$, is computationally feasible.

\section*{Acknowledgments} 
The work was supported by the Israel Science Foundation (ISF), grant ISF-589/18.
The authors would like to thank Aryeh Kontorovich and Steve Hanneke
for valuable remarks. Helpful comments of the anonymous referees are gratefully acknowledged.

\section*{Appendix}
Throughout the proofs we use various generic positive constants, not necessarily the same each
time they are used even within a single equation.

We first introduce several notations that will be used throughout the proofs. Let $|{\bf a}|_2$ be the Euclidean norm of a vector ${\bf a}$, $|A|_2$ the operator norm of a matrix $A$ and $|A|_F$ its Frobenius norm. Denote $||g||_{L_2}=(\int_\cX g^2(\bx)d\bx)^{1/2}$ for a standard $L_2$-norm of a function $g$ and $||g||_{L_2(P_X)}=(\int_\cX g^2(\bx)dP_X(\bx))^{1/2}$ for the $L_2$-norm of $g$ weighted w.r.t. the marginal distribution $P_X$ of $\bX$. In addition, the $L_\infty$-norm  $||g||_\infty=\sup_{\bx \in \cX} |g(\bx)|$.

\subsection*{Proof of Theorem \ref{th:lower}}
It is obvious that feature selection and classification in multiclass case cannot be simpler than in
binary. Formally, binary logistic classification may be viewed as a degenerate case of multinomial logistic classification with $L>2$, where  without loss of generality $p_l=0,\;l=2,\ldots,L-1$ corresponding to $\bbeta_l=-{\bf \infty},\;l=2,\ldots,L-1$ (see Section \ref{subsec:multinom}).  
Define then a subset 
$\widetilde{\cC}_L(d_0)
=\{\eta(\bx) \in \cC_L(d_0):\; \bbeta_l={\bf -\infty},\;l=2,\ldots,L-1\}$. 
Thus,
$$ 
\inf_{\widetilde \eta}\sup_{\eta^* \in \cC_L(d_0),~ P_X}
\cE(\widetilde{\eta},\eta^*) \geq \inf_{\widetilde \eta}\sup_{\eta^* \in {\widetilde \cC}_L(d_0),~ P_X}\cE(\widetilde{\eta},\eta^*) =
\inf_{\widetilde \eta}\sup_{\eta^* \in \cC_2(d_0),~ P_X}\cE(\widetilde{\eta},\eta^*)
$$
and using the results of Abramovich and Grinshtein (2019, Section 6) for binary classification we have $\inf_{\widetilde \eta}\sup_{\eta^* \in \cC_2(d_0),~ P_X}\cE(\widetilde{\eta},\eta^*) > C {\sqrt \frac{d_0 \ln \frac{de}{d_0}}{n}}$ for some $C>0$.

On the other hand, for a {\em given} model $M$ of size $d_0$, consider the corresponding set of $d_0$-dimensional linear $L$-class classifiers $\cC^M_L=\{\eta(\bx) \in \cC_L(d_0): B \in \cb_M\}$. Obviously,
$$
\inf_{\widetilde \eta}\sup_{\eta^* \in \cC_L(d_0),~ P_X}
\cE(\widetilde{\eta},\eta^*) \geq \inf_{\widetilde \eta}\sup_{\eta^* \in \cC^M_L, ~P_X}\cE(\widetilde{\eta},\eta^*).
$$
From the general results of Theorem 5 of Daniely {\em et al.}  (2015), it follows that
\be \label{eq:natbound}
\inf_{\widetilde \eta}\sup_{\eta^* \in \cC_L(d_0),~ P_X}\cE(\widetilde{\eta},\eta^*) \geq C {\sqrt \frac{d_N(\cC^M_L)}{n}}
\ee
for some $C>0$, where $d_N(\cC^M_L)$ is Natarajan dimension of
$C^M_L$.
Natarajan dimension is one of common multiclass extensions of 
VC-dimension (Natarajan, 1989) and (\ref{eq:natbound}) generalizes the corresponding well-known results for binary
classification derived in terms of VC (e.g., Devroye, Gy\"ofri and Lugosi, 1996, Chapter 14).

To complete the proof we use the bounds for Natarajan dimension of the set of $d_0$-dimensional linear $L$-class classifiers established in Daniely {\em et al.} (2012, Theorem 3.1), namely,
$d_0(L-1) \leq d_N(\cC^M_L) \leq O\left(d_0 L \ln(d_0 L)\right)$.

\subsection*{Proof of Theorem \ref{th:upper_low-noise}}
Let $KL(\bp_1,\bp_2)=\sum_{l=1}^L p_{1l}\ln\left(\frac{p_{1l}}{p_{2l}}\right)$ and $H^2(\bp_1,\bp_2)=\frac{1}{2}\sum_{l=1}^L(\sqrt{p_{1l}}-\sqrt{p_{2l}})^2$  be respectively the Kullback-Leibler divergence and the square Hellinger distance between two multinomial distributions with success probabilities
vectors $\bp_1$ and $\bp_2$. 
Let also $d_{KL}(f_{B_1},f_{B_2})=\int KL(\bp_1(\bx),\bp_2(\bx))
dP_X(\bx)$ and $d^2_H(f_{B_1},f_{B_2})=\int H^2(\bp_1(\bx),\bp_2(\bx))dP_X(\bx)$ be the corresponding Kullback-Leibler divergence and square Hellinger distance between $f_{B_1}$ and
$f_{B_2}$.

One can verify that for $\bp_1$ and $\bp_2$ satisfying Assumption (A),
$KL(\bp_1,\bp_2) \leq \frac{4(1-\delta)^2}{\delta^2} H^2(\bp_1,\bp_2)
$ and, therefore, $d_{KL}(f_{B_1},f_{B_2}) \leq \frac{4(1-\delta)^2}{\delta^2}d^2_H(f_{B_1},f_{B_2})$.  

A common approach to derive the upper bounds for misclassification risk is to convert them to the bounds of some related surrogate risk (see Section \ref{sec:intr}) which can be established by various existing methods. See, e.g., Zhang (2004ab),  Bartlett, Jordan and McAuliffe (2006), \'Avila Pires and Szepesv\'ari (2016) among many others.
 
Thus, utilizing the results of Section 5.2 of \'Avila Pires and Szepesv\'ari (2016) for multiclass logistic regression corresponding to the logistic surrogate loss and applying
then their Theorem 3.11 with the calibration function $\delta'(\epsilon)=0.5\left((1-\epsilon)\ln(1-\epsilon)+(1+\epsilon)\ln(1+\epsilon)\right) \geq 0.5\epsilon^2$ and
$\alpha'=\frac{\alpha}{\alpha+1}$ implies that under the low-noise condition (\ref{eq:low-noise})-(\ref{eq:low-noise_equiv}),
\be \label{eq:KL}
\cE(\widehat{\eta}_{\widehat M},\eta^*) \leq C 
\left(Ed_{KL}(f_B,f_{{\widehat B}_{\widehat M}}\right)^{\frac{\alpha+1}{\alpha+2}} \leq C \left(\frac{1}{\delta^2}~ E d^2_H(f_B,f_{{\widehat B}_{\widehat M}}\right)^{\frac{\alpha+1}{\alpha+2}}
\ee
and it is, therefore,  sufficient to bound the square Hellinger risk $Ed^2_H(f_B,f_{{\widehat B}_{\widehat M}})$.

We will show now that the penalty (\ref{eq:optpen}) falls within a general
class of penalties considered in Theorem \ref{thm:YB_extension} from the supplementary material which extends Theorem 1 of Yang and Barron (1998) under weaker conditions. Using this result, we find an upper bound for $Ed^2_H(f_B,f_{\widehat B_{\widehat M}})$. 

It is easy to verify that
\be \label{eq:h2lq}
H^2(\bp_1,\bp_2) \geq \frac{1}{8}|\bp_1-\bp_2|_2^2
\ee
Furthermore, using the standard inequality $\ln(1+t) \leq t$, under Assumption (A) we have
\be \label{eq:lnf}
\begin{split}
|\ln f_{B_2}(\bx,y)-\ln f_{B_1}(\bx,y)|&=\left|\sum_{l=1}^L
\xi_l \ln \frac{p_{2l}(\bx)}{p_{1l}(\bx)}\right| \leq \max_{1 \leq l \leq L} \left|\ln \frac{p_{2l}(\bx)}{p_{1l}(\bx)}\right| \\ 
& 
\leq \frac{1}{\delta}\max_{1 \leq l \leq L}|p_{2l}(\bx)-p_{1l}(\bx)|,
\end{split}
\ee
where recall that $\bxi \in \{0,1\}^L$ is the indicator vector
assigned to $y$.

For a given model $M$ consider the set of coefficient matrices $\cb_M$ defined in Section \ref{subsec:mle}.
One can easily verify that under Assumption (A), for any $B_1, B_2 \in \cb_M$ with columns $\bbeta_{1l}$'s and $\bbeta_{2l}$'s respectively and the corresponding probability vectors $\bp_1(\bx), \bp_2(\bx)$ 
\be \label{eq:equivnorm}
\delta(1-\delta)\left|(\bbeta_{2l}-\bbeta_{1l})^T \bx\right| \leq 
\left|p_{2l}(\bx)-p_{1l}(\bx)\right| \leq
\frac{1}{4}\left|(\bbeta_{2l}-\bbeta_{1l})^T \bx\right|
\ee
for all $l=1,\ldots,L-1$ and any $\bx \in \cX$.

Since $X_j$ are linearly independent, the matrix $G=E_{\bX}(\bX\bX^T)$ is positive definite. Consider the weighted Frobenius matrix norm
$|B|_G=\sqrt{tr(B^T G B)}$.
In particular, (\ref{eq:equivnorm}) implies
\be \label{eq:l2norm}
\sum_{l=1}^L \left||p_{2l}-p_{1l}\right||^2_{L_2(P_X)}  \geq \delta^2(1-\delta)^2 \sum_{l=1}^{L-1}(\bbeta_{2l}-\bbeta_{1l})^T G (\bbeta_{2l}-\bbeta_{1l}) = \delta^2 (1-\delta)^2 |B_1-B_2|^2_G
\ee
(recall that $\bbeta_{1L}=\bbeta_{2L}={\bf 0}$).

For each matrix $B_0 \in \cb_M$ consider the corresponding
Hellinger ball $\cH_{f_{B_0},r}=\{f_B: d_H(f_B,
f_{B_0}) \leq r, \;B \in \cb_M\}$.
From (\ref{eq:h2lq}) and (\ref{eq:l2norm}) it then follows that if $f_B \in \cH_{f_{B_0},r}$, the corresponding $B \in \cb_M$ lies in the
matrix ball  $\cb_{B_0,r'}=\{B \in \mathbb{R}^{|M| \times L}:
|B-B_0|_G \leq r'\}$ with $r'=\frac{2\sqrt{2} r}{\delta(1-\delta)
}$. 
 
Furthermore, for any $\bx$ and any $1 \leq l \leq L-1$, (\ref{eq:equivnorm}) and Cauchy–-Schwarz inequality imply
that 
\be \label{eq:lipschitz}
\sum_{l=1}^{L-1} \left|p_{2l}(\bx)-p_{1l}(\bx)\right|^2 \leq
\frac{1}{4} \sum_{l=1}^{L-1} (\bbeta_{2l}-\bbeta_{1l})^T G (\bbeta_{2l}-\bbeta_{1l}) \cdot |G^{-1/2} \bx |^2_2 = \frac{1}{4}|B_1-B_2|_G^2 \cdot | G^{-1/2} \bx |^2_2
\ee

Let $N(\cb_{B_0,r'},|\cdot|_G,\epsilon)$ be the $\epsilon$-covering number of
$\cb_{B_0,r^{\prime}}$ w.r.t. the $|\cdot|_G$ norm. Note that since  $\bbeta_{L} = 0$, the dimension of the vector space containing $\cb_{B_0,r^{\prime}}$ is $(L-1)|M|$. We can use then the well-known results for the covering number of a ball to have
$$
N(\cb_{B_0,r^{\prime}},|\cdot|_G,\epsilon) = N(\cb_{B_0,r^{\prime}},r^{\prime} |\cdot|_G,\epsilon / r^{\prime}) \le \left(1 + \frac{2r^{\prime}}{\epsilon}\right)^{(L-1)|M|} \le \left(\frac{3r^{\prime}}{\epsilon}\right)^{(L-1)|M|}
$$
(see, e.g., Wainwright, 2019, Example 5.8).

Consider now the bracketing number $N_{[]} (\cF_{f_{B_0},r}, \| \cdot \|_{L_2}, \epsilon)$, where $\cF_{f_{B_0},r} = \{\log f_{B}: f_{B} \in \cH_{f_{B_0},r}\}$.  By \eqref{eq:lnf} and  \eqref{eq:lipschitz}, we have
$$
|\ln f_{B_2}(\bx,y)-\ln f_{B_1}(\bx,y)| \le \frac{1}{2 \delta}|B_1-B_2|_G \cdot | G^{-1/2} \bx |_2.
$$
Let $\{B_{k},\;k=1,\ldots,N(\cb_{B_0,r'},|\cdot|_G, \delta \epsilon)\}$ be the cover set of $\cb_{B_0,r^{\prime}}$ w.r.t. the $|\cdot|_G$ norm. Define $g_{k}^{L}(\bx,y) = \log f_{B_{k}}(\bx,y) - \frac{\epsilon}{2} | G^{-1/2} \bx |_2$ and $g_{k}^{U}(\bx,y) = \log f_{B_{k}}(\bx,y) + \frac{\epsilon}{2} | G^{-1/2} \bx |_2$.
For each pair we have
$$
    \| g_{k}^{U} - g_{k}^{L} \|_{L_2} = \epsilon \| |G^{-1/2}\bx|_{2} \|_{L_{2}} = \epsilon \sqrt{E_\bX(\bX^{\top} G^{-1}\bX)} = \epsilon
$$
Finally, for any \(\log f_{B} \in \cF_{B, r^{\prime}}\), take $[g_{k}^{L}, g_{k}^{U}]$ such that $|B - B_{k}|_{G} < \delta \epsilon$. Therefore, 
\be \nonumber
\begin{split}
g_{k}^{U}(\bx,y) - \log f_{B}(\bx,y) \ge \frac{\epsilon}{2} | G^{-1/2} \bx |_2  - \frac{1}{2\delta} |B_{k} - B|_{G}  \cdot |G^{-1/2} \bx |_2\ge 0\\
g_{k}^{L}(\bx,y) - \log f_{B}(\bx,y) \le \frac{1}{2\delta} |B_{k} - B|_{G} \cdot | G^{-1/2} \bx |_2 - \frac{\epsilon}{2} | G^{-1/2} \bx |_2  \le  0,
\end{split}
\ee
which imply that $g_{k}^{L}(\bx,y) \le \log f_{B}(\bx,y) \le g_{k}^{U}(\bx,y)$. Hence, $\{[g_{k}^{L}, g_{k}^{U}]\}$ are $\epsilon$-brackets that cover $\cF_{B, r}$ under $\| \cdot \|_{L_2}$, so
\[
N_{[]} (\cF_{f_{B_0},r}, \| \cdot \|_{L_2}, \epsilon) \le N(\cb_{B_0,r'},|\cdot|_G, \delta \epsilon) \le \left(\frac{3r^{\prime}}{\delta \epsilon}\right)^{(L-1)|M|} = \left(\frac{6 \sqrt{2}}{\delta^2 (1-\delta)}~\frac{r}{\epsilon}\right)^{(L-1)|M|}
\]

The considered family of sparse multinomial logistic regression models satisfies then Assumption (D) (see supplementary material) with $A_M=\frac{18\sqrt{2}}{\delta^2 (1-\delta)}$ and $m_M=(L-1)|M|$. Note also that by Assumption (A),  \(| \sum_{l=1}^{L} \xi_{l} \ln p_l(\bx)| \le \max_{1 \leq l \leq L} | \ln p_{l}(\bx)| \le \log(1/\delta)\) for all \(\bx \in \cX\). Apply 
now Theorem \ref{thm:YB_extension}  from supplementary material for a penalized maximum likelihood model
selection procedure (\ref{eq:binpen}) with a complexity
penalty $Pen(|M|)= C_1~m_M \ln A_M+
C_2 \cdot C_M \leq \tilde{C}_1 (L-1)|M|+C_2 |M| \ln\left(\frac{de}{|M|}\right)$, where 
$C_M=|M|\ln\left(\frac{de}{|M|}\right)$. Thus,
$$
Ed^2_H(f_{\widehat B_{\widehat M}},f_B) \leq C(\delta)~ \frac{Pen(d_0)}{n} 
\leq C(\delta)~ \frac{(L-1)d_0+d_0 \ln\left(\frac{de}{d_0}\right)}{n}
$$
that together with (\ref{eq:KL}) complete the proof.

\subsection*{Proof of Theorem \ref{th:slope}} 
First, recall that from (\ref{eq:KL}) it follows that
$$
\cE(\widehat{\eta}_{gS},\eta^*) \leq C 
\left(Ed_{KL}(f_B,f_{{\widehat B}_{gS}}\right)^{\frac{\alpha+1}{\alpha+2}} 
$$
and, thus, it is sufficient to bound the Kullback-Leibler risk $Ed_{KL}(f_B,f_{{\widehat B}_{gS}})$.
For this purpose, we extend the
corresponding results of Alquier, Cottet and Lecu\'e (2019) for logistic Slope to its group analogue in multinomial logistic regression model.

As we have mentioned, the solution of (\ref{eq:groupSlope}) 
satisfies the symmetric constraint $\sum_{l=1}^L \widehat{\bbeta}_{gS,l}={\bf 0}$. 
Let $\theta_l(\bx)=\bbeta_l^T \bx,\;l=1,\ldots,L$ with 
the constraint $\sum_{l=1}^L \theta_l(\bx)=0$. Thus,
$p_l(\bx)=e^{\theta_l(\bx)}/\sum_{l'=1}^L e^{\theta_{l'}(\bx)}$
and in terms of $\btheta(\bx)$, the likelihood (\ref{eq:likelihood}) is $\ell(\btheta(\bx))=\sum_{l=1}^L y_l \theta_l(\bx)-\ln \left(\sum_{l'=1}^Le^{\theta_{l'}(\bx)}\right)$ which is Lipschitz with constant
$2$, i.e. $|\ell(\btheta_1(\bx))-\ell(\btheta_2(\bx))| \leq 2 |\btheta_1(\bx)-\btheta_2(\bx)|_2$. 
Furthermore,
similar to Lemma 1 of Abramovich and Grinshtein (2016) for binary logistic regression, re-writing the Kullback-Leibler divergence $KL(\bp_1(\bx),\bp_2(\bx))$ in terms of $\btheta(\bx)$ and expanding it in (multivariate) Taylor series, one can verify that under Assumption (A), 
$KL(\btheta_1(\bx),\btheta_2(\bx))  \geq \frac{1}{2\delta^2} |\btheta_1(\bx)-\btheta_2(\bx)|^2$ and, therefore, $d_{KL}(f_{B_1},f_{B_2}) \geq \frac{1}{2\delta^2} \sum_{l=1}^L ||\theta_{1l}(\bx)-\theta_{2l}(\bx)||^2_{L_2}$ (a multivariate analogue of Bernstein condition in terminology of Alquier, Cottet and Lecu\'e, 2019). Lipschits and Bernstein conditions allow us to apply the general approach of Alquier, Cottet and Lecu\'e (2019) and to extend their results to multinomial logistic group Lasso and group Slope. In particular, Assumption (A) corresponds to the bounded case considered there.

Let $\cB$ be a set of matrices $B$ with zero mean rows, i.e.
$\cB=\{B \in \mathbb{R}^{d \times L}: B {\bf 1}={\bf 0}\}$.
For a given regression coefficients matrix $B \in \cB$ with (zero mean) rows $B_{j\cdot}$, define its group Slope norm
$|B|_\lambda=\sum_{j=1}^d \lambda_j |B|_{(j)}$, where
recall that $|B|_{(1)} \geq \ldots \geq |B|_{(d)}$ are the
descendingly ordered $l_2$-norms of $B_{j\cdot}$'s,
and consider the corresponding unit ball $\cb_\lambda$. 
    
To derive an upper bound on \(Ed_{KL}(f_B,f_{{\widehat B}_{gS}})\) 
we define the following quantities along the lines of  Alquier, Cottet and Lecu\'e (2019).

Let $\widehat{Rad}(\cb_\lambda)$ be the empirical Rademacher complexity of $\cb_{\lambda}$, namely, 
$$
\widehat{Rad}(\cb_\lambda)=E_\Sigma\left\{\frac{1}{\sqrt n}\sup_{B 
\in \cb_\lambda} \sum_{i=1}^n \sum_{l=1}^{L} \sigma_{il} \bbeta_l^T \bX_i \Big|\bX_1=\bx_1,\ldots,\bX_n=\bx_n \right\}= E_\Sigma \left\{\frac{1}{\sqrt n}\sup_{B 
\in \cb_\lambda} tr(\Sigma B^T X^T)\right\},
$$
where the elements $\sigma_{il}$'s of $\Sigma \in \mathbb{R}^{n \times L}$ are i.i.d.
Rademacher random variables with $P(\sigma_{il}=1)=P(\sigma_{il}=-1)=1/2$, and 
$$
Rad(\cb_\lambda)=E_X\left\{\widehat{Rad}(\cb_\lambda)\right\}
$$
be the Rademacher complexity of $\cb_{\lambda}$.

Define a {\em complexity function}
$$ 
\complexity(\rho) = \sqrt{\frac{C_0 Rad(\cb_{\lambda}) \rho}{2 \delta^{2} \sqrt{n}}},\;\;\;\rho>0,
$$
where the exact value of $C_0>0$ is specified in Alquier, Cottet and Lecu\'e (2019).

Let $\mathcal{M}(\rho) = \{B \in \cB : |B|_{\lambda} = \rho,~ \sum_{l=1}^L ||B_{\cdot l}^T \bx||^2_{L_2} \le r^2(2\rho)\}$.  For a given matrix $B \in \cB$ define $\Gamma_{B}(\rho) = \bigcup_{B' : B' \in \cB,~| B' - B|_{\lambda} < \frac{\rho}{20}} \partial |\cdot |_{\lambda} (B')$, where the subdifferential $\partial |\cdot |_{\lambda} (B')=\{B'' \in \cB: |B'+B''|_\lambda - |B'|_\lambda  \geq tr((B')^T B'')\}$. The {\em sparsity parameter} is 
   $$
    \sparsity(\rho) = \inf_{B' \in \mathcal{M}} \sup_{H \in \Gamma_{B}(\rho)} <H,B'>~ = \inf_{B' \in \mathcal{M}} \sup_{H \in \Gamma_{B}(\rho)} tr(H^\top B')
    $$

Let $B \in \cB$ be $d_0$-sparse and define
\be \label{eq:rho*}
\rho^*=\frac{C_0}{800 \delta^{2}}~ \frac{Rad(\cb_{\lambda}) \left(\sum_{j=1}^{d_0} 
\lambda_j/\sqrt{j}\right)^2}{\sqrt n}
\ee
A straightforward extension of Lemma 4.3 of Lecu\'e and Mendelson (2018) for matrices implies that $\Delta(\rho^*) > \frac{4}{5}\rho^*$ and, therefore, we can apply the following
Lemma \ref{lem:ac1_extension}, which can be viewed as an extension of Theorem 2.2 (or more general Theorem 9.2)  of Alquier, Cottet and Lecu\'e (2019) for our case~:
\begin{lemma} \label{lem:ac1_extension}
Let $B \in \cB$ be  $d_0$-sparse and let $\lambda_j$'s be such that  \(Rad(\cb_{\lambda}) \leq \frac{7}{720} \sqrt{n}\).
If \(\rho^{*}\) defined in (\ref{eq:rho*}) satisfies
$\Delta(\rho^*) \geq \frac{4}{5} \rho^*$, then
    \be \label{eq:edkl}
        E d_{KL}(f_B, f_{\widehat{B}_{gS}})  \le C(\delta) \left(\sum_{j=1}^{d_0} 
\frac{\lambda_j}{\sqrt{j}}\right)^2
\ee
   for some constant $C(\delta)$ depending on $\delta$. 
\end{lemma}

To satisfy the conditions of Lemma \ref{lem:ac1_extension} and to complete the proof using (\ref{eq:lambda_slope}), we need to find an upper bound for the Rademacher complexity $Rad(\cb_{\lambda})$~:
\begin{lemma} \label{lem:max}
\be \label{eq:rad}
 Rad(\cb_\lambda) \leq C
\max_{1 \leq  j \leq d} \frac{\sqrt{L+\ln\left(\frac{d}{j}\right)}}{\lambda_j},
\ee
where the exact constant $C$ is given in the proof.
\end{lemma} 

\ignore{
Hence, combining (\ref{eq:rho*}), (\ref{eq:edkl}) and (\ref{eq:rad}) implies
\be \label{eq:edklfin}
E d_{KL}(f_B, f_{\widehat{B}_{gS}}) \leq C(\delta) \frac{1}{n}
\frac{C}{800 \delta^{2}}~ (\sum_{j=1}^{d_0} 
\lambda_j/\sqrt{j})^2 Rad(\cb_{\lambda}).
\ee
}

\subsubsection*{Proof of Lemma \ref{lem:ac1_extension}}  
The proof is an extension of the proof of Theorem 9.2 in the supplementary material of Alquier, Cottet and Lecu\'e (2019) for the multiclass framework.
    In a slightly more general version of Proposition 9.1 of Alquier, Cottet and Lecu\'e (2019)  we define the following event \(\Omega_{0}^{t}\) for $t \geq 1$:
    \be \nonumber
    \begin{split}
    \Omega_{0}^{t}= & \left\{ \forall B^{\prime} \in \cB,  \left|\frac{1}{n}(\ell(B^{\prime})-\ell(B)) -\mathbb{E}\left[(\ell(B^{\prime})-\ell(B)\right]\right| \right. \\
    & \le \left. \frac{7}{20} \delta^{2}\max\left(r\left(2\max\left(|B^{\prime} - B |_{\lambda} ,t\rho^{*}\right)\right)^{2}, 
\sum_{l=1}^L ||(B_{\cdot l}-B^{\prime}_{\cdot l})^T \bx ||^2_{L_2}\right) \right\}.
    \end{split}
     \ee
     (Proposition 9.1 of Alquier, Cottet and Lecu\'e 2019 considers only \(\Omega_{0}^{1}\)).
    
     As stated above, Assumption (A) implies the required Bernstein condition. The condition 
     \(Rad(B_{\lambda}) \leq \frac{7}{720}\sqrt{n}\) is needed for adjusting the scale of the norm w.r.t. to the loss required in Theorem 9.2 of Alquier, Cottet and Lecu\'e (2019).  Under these two conditions, we can follow the proof of Proposition 9.1 of Alquier, Cottet and Lecu\'e (2019) to get
    $$
    d_{KL}(f_B,f_{\widehat{B}_{gS}}) \le 2\delta^{2}r\left(2\rho^{*}\right)^{2} \le C \frac{\rho^* Rad(\cb_{\lambda})}{\sqrt{n}}.
    $$
on the event \(\Omega_{0}^{1}\) provided $\Delta(\rho^*) \geq \frac{4}{5}\rho^*$.
    
To extend the proof for \(t > 1\), note that \(t\rho^* \ge \rho^*\). Since \(\sparsity(\rho^*) \ge \frac{4}{5} \rho^*\), when \(|B^{\prime} - B|_{\lambda} \ge t\rho^* \ge \rho^* \) we still have \( \sparsity\left(|B^{\prime} - B|_{\lambda}\right) \ge \frac{4}{5}|B - B^{*}|_{\lambda} \) (see Lemma A.1 in Alquier, Cottet and Lecu\'e, 2019). Thus, following the arguments of Proposition 9.1,  on the event \(\Omega_{0}^{t}\) we have
    \[
    d_{KL}(f_B,f_{\widehat{B}_{gS}}) \le 2\delta^{2}r\left(2\rho^{*}\right)^{2}t \le C \frac{\rho^* Rad(\cb_{\lambda})}{\sqrt{n}} t.
    \]
     
    To bound the probability of \(\Omega_{0}^{t}\), we follow Proposition 9.3 of Alquier, Cottet and Lecu\'e (2019). We consider  the subsets \(F_{j,i}= \left\{B : \rho_{j-1} \le |B^{\prime} - B|_{\lambda} \le \rho_{j},  r_{i-1}(\rho_{j}) \le
    \sum_{l=1}^L ||(B_{\cdot l}-B^{\prime}_{\cdot l})^T \bx ||^2_{L_2}
    \le r_{i}(\rho_{j}) \right\}\), where $\rho_j=2^j \rho^*$
    and $r_i(\rho)=2^{i}r(\rho),\;i,j=0,1,\ldots$.
   
      Replace \(\rho_{j}\) with \(t\rho_{j}\) and go along the lines of the proof of Proposition 9.3 of Alquier, Cottet and Lecu\'e (2019) with the extended contraction inequality for Rademacher complexities for {\em vector}-valued Lipschitz functions of Maurer (2016) to get
    \[
    P\left(\Omega_{0}^{t}\right)\ge1-2\sum_{j=0}^{\infty}\sum_{i\in I_{j}}\exp\left(-\frac{1}{48}\widetilde{C}(\delta) \frac{7}{20}\delta^{2} n\left(2^{i}r\left(t2^{j}\rho^{*}\right)\right)^{2}\right)
    \]
    where 
    $$
    I_{j} = \{1\} \cup \left\{i \in \mathbb{N}: r_{i-1}(\rho_j) \le \min\left\{\ln\left(\frac{1-\delta}{\delta}\right), \rho_{j} \sup_{B: |B|_{\lambda} = 1} \sqrt{\sum_{l=1}^L ||B_{\cdot l}^T \bx||_{L_2}^2} \right\}\right\},
    $$
       
     Thus,
    $$
    P\left(d_{KL}(f_{\hat{B}}, f_{B})\ge2\delta^{2}r\left(2\rho^{*}\right)^{2}t\right) \le2\sum_{j=0}^{\infty}\sum_{i\in I_{j}}\exp\left(-\frac{1}{48} \widetilde{C}(\delta) \frac{7}{20}\delta^{2} n\left(2^{i}r\left(2^{j}t\rho^{*}\right)\right)^{2}\right)
    $$
    To complete the proof we use  \(\exp(- \alpha x) \le \alpha^{-1} \exp(-x),\;  x \ge 1,\; \alpha \ge 1\). 
    \ignore{
    Under the conditions of the lemma,  \(r\left(\rho\right)=\sqrt{\frac{C_0 Rad(B_{\lambda}) \rho }{2\delta^{2}\sqrt{n}}} \leq \sqrt{\frac{7C_0\rho }{1440\delta^{2}}}\).
    }
 Let \(\sqrt{n} Rad(\cb_{\lambda}) \rho^{*}>1\). For \(t > 1/\widetilde{C}(\delta)\ge 1\) we then have
    \begin{align*}
    P\left(d_{KL}(f_B,f_{\hat{B}_{gS}})\ge 2C_0\frac{Rad(B_{\lambda})\rho^{*}}{\sqrt{n}} t \right) & \le2\sum_{j=0}^{\infty}\sum_{i\in I_{j}}\exp\left(-\frac{1}{48}\widetilde{C}(\delta)\frac{7}{20}\delta^{2}\sqrt{n}2^{2i}\frac{1920}{7}\frac{1}{2\delta^{2}}2^{j}t Rad(\cb_{\lambda}) \rho^{*}\right)\\
    & = 2\sum_{j=0}^{\infty}\sum_{i\in I_{j}}\exp\left(-\widetilde{C}(\delta)\sqrt{n}2^{2i}2^{j}t Rad(\cb_{\lambda}) \rho^{*}\right) \\
     & \le2\sum_{j=0}^{\infty}\sum_{i=1}^{\infty}2^{-i}\exp\left(-\widetilde{C}(\delta)\sqrt{n}2^{j}t Rad(\cb_{\lambda}) \rho^{*}\right)\\
     & =2\sum_{j=0}^{\infty}\exp\left(-\widetilde{C}(\delta) \sqrt{n}2^{j}t Rad(\cb_{\lambda}) \rho^{*}\right) \\
     & \le 4\sum_{j=0}^{\infty}2^{-j-1}\exp\left(-\widetilde{C}(\delta)\sqrt{n} t Rad(\cb_{\lambda})\rho^{*}\right)\\
     & =4\exp\left(-\widetilde{C}(\delta)\sqrt{n}t Rad(\cb_{\lambda}) \rho^{*}\right)
    \end{align*}
Hence,    
\be \label{eq:ekl}
E d_{KL}(f_B,f_{\widehat{B}_{gS}}) \le 8\left(\frac{1}{\widetilde{C}(\delta)} + \int_{\frac{1}{\widetilde{C}(\delta)}}^{\infty} \exp(- \widetilde{C}(\delta) \sqrt{n} Rad(\cb_{\lambda}) t \rho^{*}) dt \right) C \frac{Rad(\cb_{\lambda}) \rho^*}{\sqrt{n}} \le C(\delta) \frac{Rad(\cb_{\lambda}) \rho^*}{\sqrt{n}}
\ee
Substituting $\rho^*$ from (\ref{eq:rho*}) into (\ref{eq:ekl})  under the
conditions of the lemma completes the proof.

\subsubsection*{Proof of Lemma \ref{lem:max}}
Recall that $B{\bf 1}={\bf 0}$ for $B \in \cB$.
Define the matrix  $U \in \mathbb{R}^{L \times (L-1)}$ which (orthonormal) columns are the $L-1$ eigenvectors of the
matrix $I_L-\frac{1}{L}{\bf 1}{\bf 1}^T$ corresponding to
the eigenvalue $1$. One can easily verify that $B=BUU^{\top}$.

Then,
$$
\sup_{B \in \cb_{\lambda}}tr(\Sigma B^{\top}X^{\top}) = \sup_{B \in \cb_{\lambda}}tr(X^{\top}\Sigma UU^{\top}B^{\top}) = \sup_{B \in \cb_{\lambda}}tr(K^{\top}U^{\top}B^{\top}) = \sup_{B \in \cb_{\lambda}} \sum_{j=1}^{d} K_{j \cdot}^{\top} U^{\top} B_{j\cdot},
$$
where $K=U^{\top}\Sigma^{\top}X$. Let $|K|_j=|K_{\cdot j}|_2$. 
By Cauchy-Schwartz inequality and the definition of the group
Slope norm $|B|_{\lambda}$, we have

\be \nonumber
\begin{split}
\sup_{B \in \cb_{\lambda}} \sum_{j=1}^{d} K_{j \cdot}^{\top} U B_{j\cdot} & \leq \sup_{B \in \cb_\lambda} \sum_{j=1}^d |(U B)_{j\cdot}|_2 \cdot |K_{\cdot j}|_2 = \sup_{B \in \cb_\lambda} \sum_{j=1}^d |B|_j \cdot  |K|_j \\
& =
\sup_{B \in \cb_\lambda} \sum_{j=1}^d
\lambda_j |B|_j \cdot \frac{|K|_j}{\lambda_j} 
\leq \sup_{B \in \cb_\lambda} \sum_{j=1}^d
\lambda_j |B|_{(j)} \cdot \frac{|K|_{(j)}}{\lambda_j} \\
& \leq \max_{1 \leq j \leq d} \frac{|K|_{(j)}}{\lambda_j}
\end{split}
\ee
\ignore{
and it is straightforward to construct a matrix
$B_0 \in \cb_\lambda$ that achieves the
equality in (\ref{eq:dual}). 
}
Thus, $\widehat{Rad}(\cb_\lambda)\leq E_\Sigma\left\{\max_{1 \leq j \leq d}\frac{1}{\sqrt n} \frac{|K|_{(j)}}{\lambda_j}\big|X\right\}$.

\ignore{
Let $\bx^*_j=\frac{\bx_j}{|\bx_j|_2}$ be normalized columns of $\bX$. Then,
$\frac{1}{\sqrt n}K_{\cdot j}=A_j \Sigma^{\top}\bx^*_j$, where the matrix $A_j=\frac{|\bx_j|_2}{\sqrt n}~ U^{\top}$.
}  
Let $\bx_j$ be the columns of $X$. By its definition, 
$|U|^2_F=(L-1)$ and $|U|_2=1$.  We can apply then the results of Rudelson and Vershynin (2013, p.8) to get conditionally on $\bX$
\be \nonumber
\begin{split}
P\left(\frac{|K|_j}{|\bx_j|} \geq t \sqrt{L+\ln(d/j)}\right) & = 
P\left(\frac{|K|_j}{|\bx_j|} \geq t \sqrt{L-1+\ln\left(\frac{de}{j}\right)}\right) \\
& \leq
P\left(\frac{|K|_j}{|\bx_j|} \geq \frac{t}{\sqrt 2} \sqrt{L-1}+\frac{t}{\sqrt 2} \sqrt{\ln\left(\frac{de}{j}\right)}\right) \\
& \leq P\left(\frac{|K|_j}{|\bx_j|} \geq 
\sqrt{L-1}+\frac{t}{\sqrt 2} \sqrt{\ln\left(\frac{de}{j}\right)}
\right)\\
& \leq 2 e^{-\frac{c t^2 \ln\left(\frac{de}{j}\right)}{2}}
\leq 2 \left(\frac{de}{j}\right)^{-c\frac{t^2}{2}}
\end{split}
\ee
for all $t \geq \sqrt{2}$ and a certain constant $c>0$.

Hence, by standard probabilistic arguments, for all $t \geq \max(\sqrt{2}, \frac{2}{\sqrt{c}})$
$$
P\left(\frac{|K|_{(j)}}{\lambda_j}>
t \max_{1 \leq j' \leq d}|\bx_{j'}|~ \frac{\sqrt{L+\ln(d/j')}}{\lambda_{j'}}\right) \leq 
2\binom{d}{j}\left(\frac{de}{j}\right)^{-j c\frac{t^2}{2}} \leq 2\left(\frac{de}{j}\right)^{-j \left(c\frac{t^2}{2}-1\right)} \leq 2\left(\frac{de}{j}\right)^{-jc\frac{t^2}{4}}
$$
and applying the union bound,
\be \nonumber
\begin{split}
P\left(\max_{1 \leq j \leq d}\frac{|K|_{(j)}}{\lambda_j}>
t \max_{1 \leq j' \leq d}|\bx_{j'}| \max_{1 \leq j \leq d}\frac{\sqrt{L+\ln(d/j)}}{\lambda_j}\right) & \leq
2\sum_{j=1}^d \left(\frac{de}{j}\right)^{-jc \frac{t^2}{4}} 
\leq 2\sum_{j=1}^d e^{-j c \frac{t^2}{4}} \leq 2\frac{e^{-c\frac{t^2}{4}}}
{1-e^{-c\frac{t^2}{4}}} \\
& \leq 4 e^{-c\frac{t^2}{4}}
\end{split}
\ee

Therefore, 
\be \nonumber
\begin{split}
E_\Sigma\left(\frac{\max_{1 \leq j \leq d}\frac{|K|_{(j)}}{\lambda_j}}
{\max_{1 \leq j' \leq d}|\bx_{j'}| \cdot \max_{1 \leq j \leq d}\frac{\sqrt{L+\ln(d/j)}}{\lambda_j}}\right) & =
\int_0^\infty
P\left(\max_{1 \leq j \leq d}\frac{|K|_{(j)}}{\lambda_j} > t \max_{1 \leq j' \leq d}|\bx_{j'}| \cdot
\max_{1 \leq j \leq d}\frac{\sqrt{L+\ln(d/j)}}{\lambda_j}\right)dt \\ 
& \leq \max\left(\sqrt{2}, \frac{2}{\sqrt{c}}\right) +4\int_{\max\{\sqrt{2}, \frac{2}{\sqrt{c}}\}}^{\infty} e^{- c\frac{t^2}{4}}dt 
\end{split}
\ee

Thus, $\widehat{Rad}(\cb_\lambda) \leq C \frac{1}{\sqrt n}\max_{1 \leq j' \leq d}|\bx_{j'}| \cdot \max_{1 \leq j \leq d}\frac{\sqrt{L+\ln(d/j)}}{\lambda_j}$ and by (\ref{eq:lm17}), $Rad(\cb_\lambda) \leq
C\max_{1 \leq j \leq d} \frac{\sqrt{L+\ln(d/j)}}{\lambda_j}$.



\begin{thebibliography}{99}

\bibitem{ag10}
Abramovich, F. and Grinshtein, V. (2010).
MAP model selection in Gaussian regression.
\textit{Electr. J. Statist.} \textbf{4}, 932--949.

\bibitem{ag16}
Abramovich, F. and Grinshtein, V. (2016).
Model selection and minimax estimation in generalized linear models.
\textit{IEEE Trans. Inf. Theory} \textbf{62}, 3721--3730.

\bibitem{ag19}
Abramovich, F. and Grinshtein, V. (2019).
High-dimensional classification by sparse logistic regression.
\textit{IEEE Trans. Inf. Theory} \textbf{65}, 3068--3079.

\bibitem{a13}
Agarwal, A. (2013).
Selective sampling algorithms for cost-sensitive multiclass prediction.
\textit{Proc. 30th Int. Conf. on Machine Learning,
PMLR} \textbf{28}(3), 1220--1228.

\bibitem{acl19}
Alquier, P., Cottet, V. and Lecu\'e, G. (2019).
Estimation bounds and sharp oracle inequalities of regularized procedures with Lipschitz loss functions.
\textit{Ann. Statist.} \textbf{47}, 2117-2144.

\bibitem{at07}
Audibert, J. Y. and Tsybakov, A. (2007).
Fast learning rates for plug-in classifiers.
\textit{Ann. Statist.} \textbf{35}, 608–633.


\bibitem{asg13}
\'Avila Pires, B., Szepesv\'ari, C. and Ghavamzadeh, M. (2013). 
Cost-sensitive multiclass
classification risk bounds. 
\textit{Proc. 30th Int. Conf. on Machine Learning, PMLR} \textbf{28}(3), 1391--1399.


\bibitem{as16}
\'Avila Pires, B. and Szepesv\'ari, C. (2016).
Multiclass classification calibration functions.
arxiv:1609.06385.

\bibitem{bjm06}
Bartlett, P.L., Jordan, M.I. and McAuliffe, J.D. (2006).
Convexity, classification, and risk bounds.
\textit{J. Amer. Statist. Assoc.} {\bf 101}, 138--156.

\bibitem{blt18}
Bellec, P.C., Lecu\'e, G. and Tsybakov, A. (2018).
Slope meets Lasso: improved oracle bounds and optimality.
\textit{Ann. Statist.} {\bf 46}, 3603--3642.

\bibitem{bl04}
Bickel, P. and Levina, E. (2004).
Some theory for Fisher's discriminant function, `naive Bayes', and some
alternatives where there are more variables than observations.
\textit{Bernoulli}, {\bf 10}, 989-1010.

\bibitem{brt09}
Bickel, P., Ritov, Y. and Tsybakov, A. (2009).
Simultaneous analysis of Lasso and Dantzig selector.
\textit{Ann. Statist.} {\bf 37}, 1705--1732.

\ignore{
\bibitem{b01}
Birg\'e, L. (2001).
An alternative point of view on Lepski’s method.
In: van Zwet, M.C.M.,  de Gunst C.A.J., Klaassen A.W. and van der Vaart, A. (Eds.), State of
the Art in Probability and Statistics, Festschrift for Willem R, in: Lecture Notes-Monograph Series, vol. 36, Institute of Mathematical Statistics,
113--133.
}

\bibitem{bm01}
Birg\'e, L. and Massart, P. (2001).
Gaussian model selection. \textit{J. Eur. Math. Soc.} \textbf{3},
203--268.

\bibitem{bbssc15}
Bogdan, M., van den Berg, E., Sabatti, C., Su, W. and Cand\'es, E. (2015).
SLOPE -- adaptive variable selection via convex programming.
{\em Ann. Appl. Statist.}, {\bf 9}, 1103--1140.


\bibitem{bbl05}
Boucheron, S., Bousquet, O., and Lugosi, G. (2005).
Theory of classification: a survey of some recent advances.
{\it ESAIM: Prob. Statist.} {\bf 9}, 323-375.

\bibitem{btw07}
Bunea, F., Tsybakov, A. and Wegkamp, M.H. (2007).
Aggregation for Gaussian regression.
\textit{Ann. Statist.} {\bf 4}, 1674--1697.


\bibitem{cs06}
Chen, D.-R. and Sun, T. (2006). 
Consistency of multiclass empirical risk minimization
methods based on convex loss. 
\textit{J. Mach. Learn. Res.} \textbf{7}, 2435--2447.


\bibitem{dss12}
Daniely, A., Sabato, S. and Shalev-Shwartz, S. (2012).
Multiclass learning approaches:
a theoretical comparison with implications. 
\textit{NIPS'12 Proceedings}, 485--493. 


\bibitem{dsbs15}
Daniely, A., Sabato, S., Ben-David, S. and Shalev-Shwartz, S. (2015).
Multiclass learnability and the ERM principle.
\textit{J. Mach. Learn. Res.} \textbf{16}, 2377--2404.


\bibitem{dgl95}
Devroye, L., Gy\"orfi, L. and Lugosi, G. (1996).
{\em A Probabilistic Theory of Pattern Recognition}. Springer, New York.

\bibitem{ff08} 
Fan, J. and Fan, Y. (2008). 
High-dimensional classification using feature annealed independence rules.
{\it Ann. Statist.} {\bf 36}, 2605--2637.

\ignore{
\bibitem{fht10}
Friedman, J., Hastie, H. and Tibshirani, R. (2010).
Regularization paths for generalized linear models via
coordinate descent.
\textit{J. Stat Softw.} {\bf 33},  1-–22.
}

\bibitem{km59}
Kesten, H. and Morse, N. (1959).
A property of the multinomial ditsribution.
\textit{Ann. Math. Statist.} {\bf 30}, 120--127.

\bibitem{lm17}
Lecu\'e, G. and Mendelson, S. (2017).
Sparse recovery under weak moment assumptions.
\textit{J. Eur. Math. Soc.}, {\bf 19}, 881--904.


\bibitem{lm18}
Lecu\'e, G. and Mendelson, S. (2018).
Regularization and the small-ball method I : sparse recovery. 
\textit{Ann. Statist.}, {\bf 46}, 611--641.


\bibitem{ms99}
Mammen, E. and Tsybakov, A. (1999).
Smooth discrimination analysis.
\textit{Ann. Statist.} \textbf{27}, 1808--1829.

\bibitem{m16}
Maurer, A. (2016).
A vector-contraction inequality for Rademacher complexities.
In: Ortner R., Simon H., Zilles S. (eds) \textit{Algorithmic Learning Theory. ALT 2016. Lecture Notes in Computer Science} vol {\bf 9925}, Springer, Cham, 3--17.


\bibitem{mr16}
Maximov, Yu. and Reshetova, D. (2016).
Tight risk bounds for multi-class margin classifiers.
\textit{J. Pattern Recogn. Image Analysis} \textbf{26}, 673--680.

\bibitem{mah18}
Maximov, Yu., Amini, M.R. and Harchaoui, Z. (2018).
Rademacher complexity bounds for a penalized multiclass
semi-supervised algorithm.
\textit{J. Artif. Intel. Res.} \textbf{61}, 761--786. 

\bibitem{mn89}
McCullagh, P. and Nelder, J. A. (1989).
{\em Generalized Linear Models}, 2nd ed. Chapman and Hall, London.

\bibitem{n89}
Natarajan, B.K.  (1989).
On learning sets and functions. 
\textit{Mach. Learn.} \textbf{4}, 67–-97.

\bibitem{r07}
Rigollet, P. (2007).
Generalization error bounds in semi-supervised classification under the cluster assumption.
\textit{J. Mach. Learn. Res.} \textbf{8}, 1369--1392.
 
 
\bibitem{rv13}
Rudelson, M. and Vershynin, R. (2013).
Hanson-Wright inequality and sub-gaussian concentration.
\textit{Electron. Commun. Probab.} \textbf{18} (2013), 1–-9.


\bibitem{cs15}
Su, W. and Cand\'es, E.J. (2015).
SLOPE is adaptive to unknown sparsity and asymptotically minimax. 
\textit{Ann. Statist.} {\bf 44}, 1038--1068.


\bibitem{tb07}Tewari, A. and Bartlett, P. L. (2007). 
On the consistency of multiclass classification methods.
\textit{J. Mach. Learn. Res.} \textbf{8}, 1007–-1025.

\bibitem{t04}
Tsybakov, A. (2004).
Optimal aggregation of classifiers in statistical learning.
\textit{Ann. Statist.} \textbf{32}, 135--166. 

\bibitem{vdg08}
van de Geer, S. (2008).
High-dimensional generalized linear models and the Lasso.
\textit{Ann. Statist.} \textbf{36}, 614--645.

\bibitem{v98}
Vapnik, V.N. (2000).
{\em The Nature of Statistical Learning}, 2nd ed. Springer, New York.

\ignore{
\bibitem{v18}
Vershynin, R. (2018).
{\em High-Dimensional Probability. An Introduction with Applications in Data Science}, Cambridge University Press.
}


\bibitem{vh14}
Vincent, M. and Hansen, N.R. (2014).
Sparse group lasso and high dimensional multinomial classification.
\textit{Comput. Statist. Data Anal.} \textbf{71}, 771--786.

\bibitem{w19}
Wainwright, M.J. (2019).
{\em High-Dimensional Statistics. A Non-Asymptotic Viewpoint},
Cambridge University Press.

\bibitem{yb98}
Yang, Y. and Barron, A.R. (1998)
An asymptotic property of model selection criteria.
\textit{IEEE Trans. Inf. Theory} \textbf{44}, 95--116.

\bibitem{z04a}
Zhang, T. (2004a).
Statistical behavior and consistency of classification methods based on convex risk minimization.
\textit{Ann. Statist.} \textbf{32}, 56--85.


\bibitem{z04b}
Zhang, T. (2004b). 
Statistical analysis of some multi-category large margin classification methods. 
\textit{J. Mach. Learn. Res.} \textbf{5}, 1225–-1251.


\end{thebibliography}
\end{document}


\maketitle


In supplementary material we present and prove  Theorem \ref{thm:YB_extension}, which is a stronger version of Theorem 1 of Yang and Barron (1998), where their Assumption 1 on the covering numbers of log densities classes in sup-norm is replaced by a similar condition on the bracketing numbers w.r.t. $L_2$-norm (see Assumption (D) below).

Let $f(\bx)$ be the (unknown) joint density of \(\bx = (\bx^{\prime}, \bx^{\prime\prime})\), where \(\bx^{\prime} \in \cX^{\prime} \subseteq \mathbb{R}^{d_1}, \bx^{\prime\prime} \in \cX^{\prime\prime} \subseteq \mathbb{R}^{d_2}\). Then,  \(f(\bx) = \varphi(\bx^{\prime\prime} \mid \bx^{\prime}) h(\bx^{\prime})\), where $\varphi$ is the conditional distribution of \(\bx^{\prime\prime}\) given \(\bx^{\prime}\) and \(h\) is the marginal density of \(\bx^{\prime}\).

We now try to follow the notations of Yang and Barron (1998).
Estimate $\varphi$ by considering a sequence of finite-dimensional parametric density families $\{\varphi_k
\left(\bx^{\prime\prime} \mid \bx^{\prime}, \theta^{(k)}\right), \theta^{(k)} \in \Theta_k\},\; k \in \Gamma$, where $\Gamma$ is the set of indices of the models, and view $h$ as a nuisance parameter.

Given an i.i.d. sample $\bX_1,\cdots,\bX_n$ drawn from the unknown $f$, let $\hat{\theta}^{(k)}$ be
the MLEs for $\theta^{(k)}$ under the model $k$ and select $k$ by a penalized maximum likelihood model selection criteria of the form
\be \label{eq:selection}
\hat{k} = \arg\min_{k\in\Gamma}\left\{-\sum_{i=1}^{n}\ln \varphi_{k}\left(\bX^{\prime\prime}_{i}| \bX^{\prime}, \hat{\theta}^{\left(k\right)}\right)+\lambda_{k}m_{k}+\nu C_{k}\right\}, 
\ee
where $m_k$ is the number of parameters in model $k$, $\lambda_k$ and $\nu$ are the nonnegative
tuning parameters and the weights $C_k$ satisfy $\sum_k e^{-C_k} \leq 1$.

Let $f_k(\bx, \theta^{(k)}) = \varphi_k(\bx^{\prime\prime}| \bx^{\prime}, \theta^{(k)}) h(\bx^{\prime})$.
For a given $\theta_0^{(k)} \in \Theta_k$, consider a Hellinger ball 
$\cb(\theta_0^{(k)},r)$ of densities in family $k$ of a  radius $r$ and centered at $f_{k,\theta_0^{(k)}}$, namely,
$\cb(\theta_0^{(k)},r)=\{\theta^{(k)}: \theta^{(k)} \in \Theta_k,\;d_H(f_{k, \theta^{(k)}}, f_{k, \theta_{0}^{(k)}}) < r\}$, and let \(\cF_{k}(\theta_{0}^{(k)}, r) = \{\ln f_{k, \theta^{(k)}} : d_{H}(f_{k}, f_{k, \theta_{0}^{(k)}}) < r\}\). Denote $N_{[]}(\cF_{k}(\theta_{0}^{\left(k\right)}, r), \Vert \cdot \Vert_{L_2}, \delta)$ for the $\delta$-bracketing number of $\cF_{k}(\theta_{0}^{\left(k\right)}, r)$ w.r.t. the $L_2$-norm.

\begin{assumption}[\bf D] \label{as:D}
Assume that for each \(k\in \Gamma\), there exists \(A_{k} \ge 1\) such that for any \(\theta_{0}^{(k)}\in\Theta_{k} \), \(r>0\) and \(\delta\le0.0056r\), $$
    N_{[]}(\cF_{k}(\theta_{0}^{\left(k\right)}, r), \Vert \cdot \Vert_{L_{2}}, \delta )\le\left(\frac{A_{k}r}{3\delta}\right)^{m_{k}}
$$
\end{assumption}

\begin{theorem} \label{thm:YB_extension} Assume that there exists $M>0$ such that \(||\ln \varphi||_\infty \leq M\) and \(||\ln \varphi_{k}||_\infty \leq M,\;k \in \Gamma\) for all \(\bx^{\prime} \in \cX^{\prime}\) and consider the model selection criterion (\ref{eq:selection}) with \(\lambda_{k}\ge c_1\ln A_{k}+c_2\) and \(\nu\ge c_3\), where \(c_1,c_2,c_3\) are certain universal constants. Then, for any $h$, under Assumption (D), 
\[
\mathbb{E}_{\bX}d_{H}^{2}(f,f_{\hat{k},\hat{\theta}^{\left(\hat{k}\right)}})\le C(M)~R_{n}\left(f\right),
\]
where 
\[
R_{n}\left(f\right)=\inf_{k\in\Gamma}\left\{ \inf_{\theta^{\left(k\right)}\in\Theta_{k}} d_{KL}(f, f_{k,\theta^{\left(k\right)}})+\frac{\lambda_{k}m_{k}}{n}+\frac{\nu C_{k}}{n}\right)
\]
\end{theorem}

The proof of Theorem \ref{thm:YB_extension}
essentially repeats  the proof of Theorem 1 in Yang and Barron (1998) with their Lemma 0 been replaced by the following Lemma \ref{lem:YB_main_lemma} and, as a result, with somewhat different numerical constants:

 \begin{lemma} \label{lem:YB_main_lemma}
 Let \(A>0\), \(m\ge1\), and \(\rho>0\) and \(\rho\le A\) such that for any \(r>0\) and \(\delta\le\rho r\), 
 \[
 N_{[]}\left(\mathcal{F}(\theta, r),\left\Vert \cdot\right\Vert _{L_{2}}, \delta\right)\le\left(\frac{Ar}{\delta}\right)^{m}.
 \]
 Assume that there exists $M>0$ such that \(||\ln \varphi||_\infty \leq M\) and \(||\ln \varphi_{k}||_\infty \leq M,\;k \in \Gamma\) for all \(\bx^{\prime} \in \cX^{\prime}\).
Assume also that   \(\rho\ge\gamma/6\sqrt{10}\sqrt{1-4\gamma}\) for some \(0 \le \gamma \le \frac{1}{4}\). 

If \(\frac{\xi}{m}\ge\frac{4}{1-4\gamma}\ln\left(\frac{12 A\sqrt{10}\sqrt{1-4\gamma}}{\gamma}\right)\), then
 \[
 P\left(\exists \theta\in\Theta, \frac{1}{n}\sum_{t=1}^{n}\ln\frac{f\left(\bX_{t},\theta\right)}{f\left(\bX_{t}\right)}\ge-\gamma d_{H}^{2}\left(f,f_{\theta}\right)+\frac{\xi}{n}\right)\le C(\rho) \exp\left(-\frac{\left(1-4\gamma\right)\xi}{8}\right),
 \]
 for some \(C(\rho) > 0\).
 \end{lemma}

\noindent
\newline
{\bf Proof of Lemma \ref{lem:YB_main_lemma}}
\begin{proof}
The proof of Lemma \ref{lem:YB_main_lemma} is a modification of the proof of Lemma 0 of Yang and Barron (1998) utilizing the bracket covering similar to Theorems 5.11 and 8.13 in van de Geer (2000).

Divide the parameter space into rings \(r_{i}=2^{i/2}\sqrt{\frac{\xi}{n}}\) with  \(r_{-1}=0\), and consider \[
\Theta_{i}=\left\{ \theta\in\Theta:r_{i-1}\le d_{H}^{2}\left(f,f_{\theta}\right)\le r_{i}\right\} .
\]
If we show that for some constant \(C_1(\rho)\),
\be \label{eq:YB_ring_bound}
P\left(\exists \theta\in\Theta_{i}, \frac{1}{n}\sum_{t=1}^{n}\ln\frac{f\left(\bX_{t},\theta\right)}{f\left(\bX_{t}\right)}\ge-\gamma d_{H}^{2}\left(f,f_{\theta}\right)+\frac{\xi}{n}\right)\le C_1(\rho) \exp\left(-\frac{\left(i+1\right)\left(1-4\gamma\right)\xi}{8}\right),
\ee
then 
\[
\begin{split}
P\left(\exists \theta\in\Theta, \frac{1}{n}\sum_{t=1}^{n}\ln\frac{f\left(\bX_{t},\theta\right)}{f\left(\bX_{t}\right)}\ge-\gamma d_{H}^{2}\left(f,f_{\theta}\right)+\frac{\xi}{n}\right)	\le & \sum_{i=0}^{\infty}P\left(\exists \theta\in\Theta_{i}, \frac{1}{n}\sum_{t=1}^{n}\ln\frac{f\left(\bX_{t},\theta\right)}{f\left(\bX_{t}\right)}\ge-\gamma d_{H}^{2}\left(f,f_{\theta}\right)+\frac{\xi}{n}\right) \\
	\le & C_1(\rho) \sum_{i=0}^{\infty}\exp\left(-\frac{\left(i+1\right)\left(1-4\gamma\right)\xi}{8}\right) \\
	= & C_1(\rho) \left(\frac{1}{1-\exp\left(-\frac{\left(1-4\gamma\right)\xi}{8}\right)}\right)\exp\left(-\frac{\left(1-4\gamma\right)\xi}{8 }\right) 
\end{split}
\]

Since the condition of the lemma guarantees \(\frac{\xi}{m}\ge\frac{4 }{1-4\gamma}\ln\frac{2A}{\rho},\;\rho\le A\), we have \(\frac{\ln2}{2}\le\frac{\ln\left(2\right)}{2}m\le\frac{\left(i+1\right)\left(1-4\gamma\right)\xi}{8}.\) 
Thus, for some constant \(C(\rho)\),
\[
P\left(\exists \theta\in\Theta, \frac{1}{n}\sum_{t=1}^{n}\ln\frac{f\left(\bX_{t},\theta\right)}{f\left(\bX_{t}\right)}\ge-\gamma d_{H}^{2}\left(f,f_{\theta}\right)+\frac{\xi}{n}\right) \le C(\rho) \exp\left(-\frac{\left(1-4\gamma\right)\xi}{8}\right).
\]

We now turn to establishing \eqref{eq:YB_ring_bound}. Denote \( g_{\theta}\left(\bX_{t}\right)=\ln \varphi(\bX^{\prime\prime}_t| \bX^{\prime}_t, \theta) -\mathbb{E}_{\bX}\left(\varphi(\bX^{\prime\prime}| \bX^{\prime}, \theta)\right)\) and \(g\left(\bX_{t}\right)= \ln \varphi(\bX^{\prime\prime}_t| \bX^{\prime}_t) -\mathbb{E}_{\bX}\left(\varphi(\bX^{\prime\prime}| \bX^{\prime})\right)\). Then, \[
\frac{1}{n}\sum_{t=1}^{n}\ln\frac{f\left(\bX_{t},\theta\right)}{f\left(\bX_{t}\right)}=\frac{1}{n}\sum_{t=1}^{n}g_{\theta}\left(\bX_{t}\right)-g\left(\bX_{t}\right)+\mathbb{E}_{\bX}\left(\frac{\ln\left(f\left(\bX,\theta\right)\right)}{\ln\left(f\left(\bX\right)\right)}\right).
\]
Denote \(\mathcal{G}(\theta, r_i) = \{\ln \varphi(\cdot, \theta) - \mathbb{E}_{\bX}\left(\ln \varphi(\bX^{\prime\prime}| \bX^{\prime}, \theta)\right): d_H(f,f_\theta) < r\}\). Shifting each bracket \([f^{L}, f^{U}]\) to \([f^{L} - \ln h - \mathbb{E}_{\bX}(f^{L} - \ln h), f^{U} - \ln h - \mathbb{E}_{\bX}(f^{U} - \ln h)]\) implies 
\[
N_{[]}(\mathcal{G}(\theta, r_i), \Vert \cdot \Vert_{L_{2}}, \delta ) \le N_{[]}(\mathcal{F}(\theta, r_i), \Vert \cdot \Vert_{L_{2}}, \delta ). 
\] 
Let \(\delta_{j}\to 0\) be a monotonically decreasing sequence to be specified later such that \(\delta_{0}\le\rho r_{0}\). Let \(\left\{ \left[\tilde{g}_{s}^{j,U},\tilde{g}_{s}^{j,L}\right]\right\} _{s=1}^{N_{[]}\left(\mathcal{G}(\theta, r_i),\left\Vert \cdot\right\Vert, \delta_{j} \right)}\) be the bracketing set for \(\mathcal{G}(\theta, r_i)\) with \(\delta_{j}\). By definition, for all \(\theta\in\Theta\) and \(\delta_{j}\), 
\[
\begin{split}
& \tilde{g}_{s\left(\theta\right)}^{j,L}\left(\bx\right) \le	g_{\theta}\left(\bx\right)\le\tilde{g}_{s\left(\theta\right)}^{j,U}\left(\bx\right) \\
& \left\Vert \tilde{g}_{s\left(\theta\right)}^{j,U}-\tilde{g}_{s\left(\theta\right)}^{j,L}\right\Vert_{L_2} \le \delta_{j}.
\end{split}
\]
Let \(g_{s\left(\theta\right)}^{j,L}=\max_{r\le j}\tilde{g}_{s\left(\theta\right)}^{r,L}\) and \(g_{s\left(\theta\right)}^{j,U}=\min_{r\le j}\tilde{g}_{s\left(\theta\right)}^{r,U}\). Denote \(\Delta_{j,s\left(\theta\right)}=g_{s\left(\theta\right)}^{j,U}-g_{s\left(\theta\right)}^{j,L}\) and note that \[
\Delta_{j,s\left(\theta\right)} \ge g_{\left(s\right)}-g_{s\left(\theta\right)}^{j,L}\ge 0,
\]
where \(\left\Vert \Delta_{i,j,s\left(\theta\right)}\right\Vert \le\delta_{j}\). For all \(\bx\),
\[
0\le g_{\theta}\left(\bx\right)-\left(g_{s\left(\theta\right)}^{0,L}+\sum_{j=1}^{m}\left(g_{s\left(\theta\right)}^{j,L}\left(\bx\right)-g_{s\left(\theta\right)}^{j-1,L}\left(\bx\right)\right)\right)=g_{\theta}\left(\bx\right)-g_{s\left(\theta\right)}^{m,L}\left(\bx\right)\le\Delta_{m,s\left(\theta\right)}.
\]
Thus,
\[
\left\Vert g_{\theta}-\left(g_{s\left(\theta\right)}^{0,L}+\sum_{j=1}^{m}\left(g_{s\left(\theta\right)}^{j,L}-g_{s\left(\theta\right)}^{j-1,L}\right)\right)\right\Vert _{L_{2}}=\left\Vert g_{\theta}-g_{s\left(\theta\right)}^{m,L}\right\Vert _{L_{2}}\le\delta_{m}
\]
Consider exponentially decreasing sequence \(\delta_{j}\le a\exp\left(-bj\right)\). By Chebyshev inequality,
\[
P\left(\left|g_{\theta}\left(\bX\right)-\left(g_{s\left(\theta\right)}^{0,L}\left(\bX\right)+\sum_{j=1}^{m}\left(g_{s\left(\theta\right)}^{j,L}\left(\bX\right)-g_{s\left(\theta\right)}^{j-1,L}\left(\bX\right)\right)\right)\right|\ge\frac{1}{m}\right)\le a^{2}\exp\left(-2bm\right)m^{2}\underset{m\to\infty}{\to}0
\]
Therefore, \(g_{s\left(\theta\right)}^{0,L}+\sum_{j=1}^{\infty}\left(g_{s\left(\theta\right)}^{j,L}-g_{s\left(\theta\right)}^{j-1,L}\right)\) convergence almost surely and
\be \label{eq:g_telescope}
g_{\theta}\left(\bx\right)	=g_{s\left(\theta\right)}^{0,L}\left(\bx\right)+\sum_{j=1}^{\infty}\left(g_{s\left(\theta\right)}^{j,L}\left(\bx\right)-g_{s\left(\theta\right)}^{j-1,L}\left(\bx\right)\right).
\ee
To perform adaptive truncation, let \(K_{s}\left(\bx\right)\) be the first index \(j\) such that \(\Delta_{j,s}\left(\bx\right)>A_{j}\) (or infinity if no such index exists). Using the the truncation levels \(K_{s}\), rewrite  \eqref{eq:g_telescope} as 
\[
g_{\theta}=g_{s\left(\theta\right)}^{0,L}+\sum_{j=1}^{\infty}\left(g_{s\left(\theta\right)}^{j,L}-g_{s\left(\theta\right)}^{j-1,L}\right)1\left(K_{s\left(\theta\right)}\ge j\right)+\sum_{j=1}^{\infty}\left(g_{s\left(\theta\right)}^{j,L}-g_{s\left(\theta\right)}^{j-1,L}\right)1\left(K_{s\left(\theta\right)}<j\right).
\]
Note that for every \(\bx\),
\[
\begin{split}
\sum_{j=1}^{\infty}\left(g_{s\left(\theta\right)}^{j,L}\left(\bx\right)-g_{s\left(\theta\right)}^{j-1,L}\left(\bx\right)\right)1\left(K_{s\left(\theta\right)}\left(\bx\right)<j\right)	&=\sum_{j=K_{s\left(\theta\right)}\left(\bx\right)+1}^{\infty}\left(g_{s\left(\theta\right)}^{j,L}\left(\bx\right)-g_{s\left(\theta\right)}^{j-1,L}\left(\bx\right)\right) \\
	&=g_{\theta}\left(\bx\right)-g_{s\left(\theta\right)}^{K_{s\left(\theta\right)}\left(\bx\right),L}\left(\bx\right) \\
	&=\sum_{j=0}^{\infty}\left(g_{\theta}\left(\bx\right)-g_{s\left(\theta\right)}^{j,L}\left(\bx\right)\right)1\left(j=K_{s\left(\theta\right)}\left(\bx\right)\right).
\end{split}
\]
Thus, 
\[
\begin{split}
g_{\theta}\left(\bx\right)=&\tilde{g}_{s\left(\theta\right)}^{0,L}+\sum_{j=1}^{\infty}\left(g_{s\left(\theta\right)}^{j,L}-g_{s\left(\theta\right)}^{j-1,L}\right)1\left(K_{s\left(\theta\right)}(\bx)\ge j\right) \\
&+\sum_{j=1}^{\infty}\left(g_{\theta}\left(\bx\right)-g_{s\left(\theta\right)}^{j,L}\left(\bx\right)\right)1\left(j=K_{s\left(\theta\right)}\left(\bx\right)\right)\\ &+\left(g_{\theta}\left(\bx\right)-g_{s\left(\theta\right)}^{0,L}\left(\bx\right)\right)1\left(K_{s\left(\theta\right) = 0}\left(\bx\right)\right).
\end{split}
\]
For brevity, we write \(K_{s\left(\theta\right),t}\) for the random variable \(K_{s\left(\theta\right)}(\bX_t)\). Therefore, for \(\alpha = -\gamma d_{H}^2(f,f_{\theta}) + \frac{\xi}{n}\),
\[
\begin{split}
P\left(\exists \theta\in\Theta_{i}, \frac{1}{n}\sum_{t=1}^{n}\ln\frac{f\left(\bX_{t},\theta\right)}{f\left(\bX_{t}\right)}\ge \alpha \right) &\le	P\left(\exists \theta\in\Theta_{i}, \frac{1}{n}\sum_{t=1}^{n}g_{\theta}\left(\bX_{t}\right)-g\left(\bX_{t}\right)+\mathbb{E}_{\bX}\left(\frac{\ln\left(f\left(\bX,\theta\right)\right)}{\ln\left(f\left(\bX\right)\right)}\right)\ge-\gamma r_{i}^{2}+\frac{\xi}{n}\right) \\
& \le	P\left(\exists \theta\in\Theta_{i}, \frac{1}{n}\sum_{t=1}^{n}g_{s\left(\theta\right)}^{0,L}\left(\bX_{t}\right)-g\left(\bX_{t}\right)+\mathbb{E}_{\bX}\left(\frac{\ln\left(f\left(\bX,\theta\right)\right)}{\ln\left(f\left(\bX\right)\right)}\right)\ge-2\gamma r_{i}^{2}+\frac{\xi}{n}\right) \\
	&+P\left(\exists \theta\in\Theta_{i}, \frac{1}{n}\sum_{t=1}^{n}\sum_{j=1}^{\infty}\left(g_{s\left(\theta\right)}^{j,L}\left(\bX_{t}\right)-g_{s\left(\theta\right)}^{j-1,L}\left(\bX_{t}\right)\right)1\left\{K_{s\left(\theta\right),t}\ge j\right\}\ge\frac{1}{3}\gamma r_{i}^{2}\right) \\
	&+P\left(\exists \theta\in\Theta_{i}, \frac{1}{n}\sum_{t=1}^{n}\sum_{j=1}^{\infty}\left(g_{\theta}\left(\bX_{t}\right)-g_{s\left(\theta\right)}^{j,L}\left(\bX_{t}\right)\right)1\left\{K_{s\left(\theta\right),t}=j\right\}\ge\frac{1}{3}\gamma r_{i}^{2}\right)\\
	&+P\left(\exists \theta\in\Theta_{i}, \frac{1}{n}\sum_{t=1}^{n} \left(g_{\theta}\left(\bX_{t}\right)-g_{s\left(\theta\right)}^{0,L}\left(\bX_{t}\right)\right)1\left\{K_{s\left(\theta\right),t}=0\right\}\ge\frac{1}{3}\gamma r_{i}^{2}\right)\\
&=	P_{I}+P_{II}+P_{III}+P_{IV}.
\end{split}
\]

\begin{enumerate}
\item \textbf{Bounding \(P_{I}\)}
    
Let \(1\le s\le N_{[]}\left(\delta_{0},\mathcal{G}(\theta, r_i), \left\Vert \cdot\right\Vert \right)\). For each bracket \(\left[g_{s}^{j,U},g_{s}^{j,L}\right]\) consider 
\[
\Theta_{0}\left(s\right)=\left\{ \theta:g_{s}^{j,L}\le g_{\theta}\le g_{s}^{j,U},\theta\in\Theta_{i}\right\},
\]
and fix \(\theta_{s}\in\Theta_{0}\left(s\right)\). Let \(\theta_{0}\left(s\right)\in\Theta_{s}\) be such that for \(\varepsilon>0\),
\[
\mathbb{E}_{\bX}\left(\ln\frac{f\left(\bX,\theta_{s}\right)}{f\left(\bX,\theta_{0}\left(s\right)\right)}\right)\le\inf_{\theta\in\Theta_{0}\left(s\right)}\mathbb{E}_{\bX}\left(\ln\frac{f\left(\bX,\theta_{s}\right)}{f\left(\bX,\theta\right)}\right)+\varepsilon
\]
Therefore, for all \(\theta\in\Theta_{0}\left(s\right)\), \[
\begin{split}
\mathbb{E}_{\bX}\left(\ln\frac{f\left(\bX,\theta_{0}\left(s\right)\right)}{f\left(\bX,\theta\right)}\right)	&= \mathbb{E}_{\bX}\left(\ln\frac{f\left(\bX,\theta_{0}\left(s\right)\right)}{f\left(\bX,\theta_{s}\right)}\right)+\mathbb{E}_{\bX}\left(\ln\frac{f\left(\bX,\theta_{s}\right)}{f\left(\bX,\theta\right)}\right) \\
	&\ge\mathbb{E}_{\bX}\left(\ln\frac{f\left(\bX,\theta_{s}\right)}{f\left(\bX,\theta\right)}\right)-\inf_{\theta\in\Theta_{0}\left(s\right)}\mathbb{E}_{\bX}\left(\ln\frac{f\left(\bX,\theta_{s}\right)}{f\left(\bX,\theta\right)}\right)-\varepsilon \\
	& \ge-\varepsilon
\end{split}
\]
Note that for all \(s\), \(g_{s}^{0,L}\left(\bX\right)\le\ln\left(f\left(\bX,\theta_{0}\left(s\right)\right)\right)-\mathbb{E}_{\bX}\left(\ln\left(f\left(\bX,\theta_{0}\left(s\right)\right)\right)\right)\le g_{s\left(\theta\right)}^{0,U}\). Thus,
\[
\begin{split}
g_{s\left(\theta\right)}^{0,L}\left(\bX\right)-g\left(\bX\right)+\mathbb{E}_{\bX}\left(\frac{\ln\left(f\left(\bX,\theta\right)\right)}{\ln\left(f\left(\bX\right)\right)}\right)	&\le\frac{\ln\left(f\left(\bX,\theta_{0}\left(s\right)\right)\right)}{\ln\left(f\left(\bX\right)\right)}-\mathbb{E}_{\bX}\left(\frac{\ln\left(f\left(\bX,\theta_{0}\left(s\right)\right)\right)}{\ln\left(f\left(\bX,\theta\right)\right)}\right) \\
	&\le\frac{\ln\left(f\left(\bX,\theta_{0}\left(s\right)\right)\right)}{\ln\left(f\left(\bX\right)\right)}+\varepsilon
\end{split}
\]
and, therefore,
\[
\begin{split}
P_{I} &= P\left(\exists \theta\in\Theta_{i}, \frac{1}{n}\sum_{t=1}^{n}g_{s\left(\theta\right)}^{0,L}\left(\bX_{t}\right)-g\left(\bX_{t}\right)+\mathbb{E}_{\bX}\left(\frac{\ln\left(f\left(\bX,\theta\right)\right)}{\ln\left(f\left(\bX\right)\right)}\right)\ge-2\gamma r_{i}^{2}+\frac{\xi}{n}\right) \\
	&\le P\left(\sup_{s\le N\left(\mathcal{G}(\theta, r_i),\left\Vert \cdot\right\Vert_{L_2}, \delta_{0} \right)}\frac{1}{n}\sum_{t=1}^{n}\frac{\ln\left(f\left(\bX_{t},\theta_{0}\left(s\right)\right)\right)}{\ln\left(f\left(\bX_{t}\right)\right)}+\varepsilon\ge-2\gamma r_{i}^{2}+\frac{\xi}{n}\right) \\
	&\le\sum_{s=1}^{N_{[]}\left(\mathcal{G}(\theta, r_i),\left\Vert \cdot\right\Vert_{L_2} , \delta_{0}\right)}P\left(\frac{1}{n}\sum_{t=1}^{n}\frac{\ln\left(f\left(\bX_{t},\theta_{0}\left(s\right)\right)\right)}{\ln\left(f\left(\bX_{t}\right)\right)}\ge-2\gamma r_{i}^{2}+\frac{\xi}{n}-\varepsilon\right) \\
	&\le\sum_{s=1}^{N_{[]}\left(\mathcal{G}(\theta, r_i),\left\Vert \cdot\right\Vert_{L_2}, \delta_{0} \right)}\exp\left(-\frac{n}{2}d_{H}^{2}\left(f,f_{\theta_{0}\left(s\right)}\right)-2\gamma r_{i}^{2}+\frac{\xi}{n}-\varepsilon\right) \\
	&\le\left(\frac{Ar_{i}}{\delta_{0}}\right)^{n}\exp\left(-\frac{n}{2}\left(r_{i-1}^{2}-2\gamma r_{i}^{2}+\frac{\xi}{n}\right)-\varepsilon\right),
\end{split}
\]
where the third inequality is derived from the exponential inequality in Yang and Barron (1998). Since this holds for any \(\varepsilon>0\), we have 
\[
\begin{split}
P_{I} &=P\left(\exists \theta\in\Theta_{i}, \sum_{t=1}^{n}g_{\theta}\left(\bX_{t}\right)-g\left(\bX_{t}\right)+\mathbb{E}_{\bX}\left(\frac{\ln\left(f\left(\bX,\theta\right)\right)}{\ln\left(f\left(\bX\right)\right)}\right) \ge-2\gamma r_{i}^{2}+\frac{\xi}{n}\right) \\
&\le\left(\frac{Ar_{i}}{\delta_{0}}\right)^{n}\exp\left(-\frac{n}{2}\left(r_{i-1}^{2}-2\gamma r_{i}^{2}+\frac{\xi}{n}\right)\right).
\end{split}
\]
Finally, since \(r_{i}=2^{i/2}\sqrt{\frac{\xi}{n}}\), for \(i\ge1\), 
\[
\begin{split}
r_{i-1}^{2}-2\gamma r_{i}^{2}+\frac{\xi}{n}	&=2^{i-1}\frac{\xi}{n}-2^{i+1}\gamma\frac{\xi}{n}+\frac{\xi}{n}=\left(2^{i-1}-2^{i+1}\gamma+1\right)\frac{\xi}{n} \\
	& \ge\left(i+1\right)\left(1-4\gamma\right)\frac{\xi}{n}
\end{split}
\]
and for \(i=0\), 
\[
r_{-1}^{2}-2\gamma r_{0}^{2}+\frac{\xi}{n}=\frac{\xi}{n}-2\gamma\frac{\xi}{n}\ge\left(1-4\gamma\right)\frac{\xi}{n}.
\]
Thus,
\[
P_{I}\le\left(\frac{Ar_{i}}{\delta_{0}}\right)^{n}\exp\left(-\frac{\left(i+1\right)\left(1-4\gamma\right)\xi}{2}\right)
\]

\item \textbf{Bounding \(P_{II}\)}

Note that when \(K_{s\left(\theta\right)}\left(\bx\right)\ge j\), since \(g_{s\left(\theta\right)}^{j,L}\left(\bx\right)\) is non-decreasing in \(j\), 
\[
\begin{split}
\left|g_{s\left(\theta\right)}^{j,L}\left(\bx\right)-g_{s\left(\theta\right)}^{j-1,L}\left(\bx\right)\right|	&\le\left|g_{\theta}\left(\bx\right)-g_{s\left(\theta\right)}^{j-1,L}\left(\bx\right)\right| \\
	&\le\Delta_{j-1,s}\left(\bx\right) \\
	&\le A_{j-1}.
\end{split}
\] 
On the other hand, 
\[
    \left\Vert g_{s\left(\theta\right)}^{j,L}\left(\bx\right)-g_{s\left(\theta\right)}^{j-1,L}\left(\bx\right)\right\Vert _{L_{2}}\le\left\Vert g_{s\left(\theta\right)}^{j-1,U}\left(\bx\right)-g_{s\left(\theta\right)}^{j-1,L}\left(\bx\right)\right\Vert _{L_{2}}\le\delta_{j-1}.
\]
Thus, by Bernstein inequality, 
\[
P\left(\frac{1}{n}\sum_{t=1}^{n}\left(g_{s\left(\theta\right)}^{j,L}\left(\bX_{t}\right)-g_{s\left(\theta\right)}^{j-1,L}\left(\bX_{t}\right)\right)1\left\{K_{s\left(\theta\right)}\left(\bX_{t}\right)\ge j\right\}\ge\eta_{j}\right)\le\exp\left(-\frac{n\eta_{j}^{2}}{2\left(\delta_{j-1}^{2}+\frac{1}{3}A_{j-1}\eta_{j}\right)}\right).
\]
Let \(\eta_{j}\) such that \(\sum_{j=0}^{\infty}\eta_{j}\le\frac{1}{3}\gamma r_{i}^{2}\). Then, for \(A_{j}=9\frac{\delta_{j}^{2}}{\eta_{j+1}}\),

\[
\begin{split}
P_{II} &= P\left(\exists \theta\in\Theta_{i}, \sum_{j=1}^{\infty}\frac{1}{n}\sum_{t=1}^{n}\left(g_{s\left(\theta\right)}^{j,L}\left(\bX_{t}\right)-g_{s\left(\theta\right)}^{j-1,L}\left(\bX_{t}\right)\right)1\left\{K_{s\left(\theta\right)}\left(\bX_{t}\right)\ge j\right\}\ge \frac{1}{3} \gamma r_i^{2} \right) \\
	&\le\sum_{j=1}^{\infty}P\left(\exists \theta\in\Theta_{i}, \frac{1}{n}\sum_{t=1}^{n}\left(g_{s\left(\theta\right)}^{j,L}\left(\bX_{t}\right)-g_{s\left(\theta\right)}^{j-1,L}\left(\bX_{t}\right)\right)1\left\{K_{s\left(\theta\right)}\left(\bX_{t}\right)\ge j\right\}\ge\eta_{j}\right) \\
	& \le\sum_{j=1}^{\infty}N_{[]}\left(\mathcal{G}(\theta, r_{i}),\left\Vert \cdot\right\Vert_{L_2}, \delta_{j} \right)N_{[]}\left(\mathcal{G}(\theta, r_{i}),\left\Vert \cdot\right\Vert_{L_2}, \delta_{j-1} \right)\exp\left(-\frac{n\eta_{j}^{2}}{2\left(\delta_{j-1}^{2}+\frac{1}{3}A_{j-1}\eta_{j}\right)}\right) \\
	& \le\sum_{j=1}^{\infty}\left(\frac{Ar_{i}}{\delta_{j-1}}\right)^{m}\left(\frac{Ar_{i}}{\delta_{j}}\right)^{m}\exp\left(-\frac{n\eta_{j}^{2}}{8\delta_{j-1}^{2}}\right) \\
\end{split}
\]

\item \textbf{Bounding \(P_{III}\)}

Note that for \(K_{s\left(\theta\right)}\left(\bX_{t}\right)=j,\; j \ge 1\), \[
\left|g_{\theta}\left(\bx\right)-g_{s\left(\theta\right)}^{j,L}\left(\bx\right)\right|\le\left|g_{\theta}\left(\bx\right)-g_{s\left(\theta\right)}^{j-1,L}\left(\bx\right)\right|\le\Delta_{j-1,s}\left(\bx\right)\le A_{j-1}
\]
and 
\[
\left\Vert g_{\theta}\left(\bx\right)-g_{s\left(\theta\right)}^{j,L}\left(\bx\right)\right\Vert _{L_{2}}\le\left\Vert g_{s\left(\theta\right)}^{j,U}\left(\bx\right)-g_{s\left(\theta\right)}^{j,L}\left(\bx\right)\right\Vert _{L_{2}}\le\delta_{j}\le\delta_{j-1}
\]
Thus, by Bernstein inequality, \[
P\left(\frac{1}{n}\sum_{t=1}^{n}\left(g_{\theta}\left(\bX_{t}\right)-g_{s\left(\theta\right)}^{j,L}\left(\bX_{t}\right)\right)1\left\{K_{s\left(\theta\right)}\left(\bX_{t}\right)=j\right\}\ge\eta_{j}\right)\le\exp\left(\frac{n\eta_{j}^{2}}{2\left(\delta_{j-1}^{2}+\frac{1}{3}A_{j-1}\eta_{j}\right)}\right)
\]

Consider \(\eta_{j}\) from the bound of \(P_{II}\). Thus, for the same  \(A_{j}=9\frac{\delta_{j}^{2}}{\eta_{j+1}}\),
\[
\begin{split}
P_{III} &= P\left(\exists \theta\in\Theta_{i}, \sum_{j=0}^{\infty}\frac{1}{n}\sum_{t=1}^{n}\left(g_{\theta}\left(\bX_{t}\right)-g_{s\left(\theta\right)}^{j,L}\left(\bX_{t}\right)\right)1\left\{K_{s\left(\theta\right)}\left(\bX_{t}\right)=j\right\}\ge\frac{1}{2}\gamma r_{i}^{2}\right) \\
	&\le\sum_{j=0}^{\infty}P\left(\exists \theta\in\Theta_{i}, \frac{1}{n}\sum_{t=1}^{n}\left(g_{s\left(\theta\right)}^{j,L}\left(\bX_{t}\right)-g_{s\left(\theta\right)}^{j-1,L}\left(\bX_{t}\right)\right)1\left\{K_{s\left(\theta\right)}\left(\bX_{t}\right)\ge j\right\}\ge\eta_{j}\right) \\
	& \le\sum_{j=0}^{\infty}N_{[]}\left(\mathcal{G}(\theta, r_{i}),\left\Vert \cdot\right\Vert_{L_2}, \delta_{j} \right)N_{[]}\left(\mathcal{G}(\theta, r_{i}),\left\Vert \cdot\right\Vert_{L_2}, \delta_{j-1} \right)\exp\left(\frac{n\eta_{j}^{2}}{2\left(\delta_{j-1}^{2}+\frac{1}{3}A_{j-1}\eta_{j}\right)}\right) \\
	& \le\sum_{j=1}^{\infty}\left(\frac{Ar_{i}}{\delta_{j-1}}\right)^{m}\left(\frac{Ar_{i}}{\delta_{j}}\right)^{m}\exp\left(-\frac{n\eta_{j}^{2}}{8\delta_{j-1}^{2}}\right)
\end{split}
\]

\item \textbf{Bounding \(P_{IV}\)}

Let \(\eta_{0} \le \min\left\{1/3 \gamma r_i, M \sqrt{2 \ln (2)\frac{mi}{n}+ \frac{(i+1)(1-4\gamma)\xi}{n}}\right\}\). Since \(\sup_{\bx} | \ln f(\bx) - \ln h(\bx^{\prime}) | = |\ln \phi(\bx^{\prime\prime}|\bx^{\prime})| \le M\), we have
    \[
    0 \le \sup_{\bx} g(\bx) - g_{s(\theta)}^{0,L}(\bx) \le \sup_{\bx} | g(\bx) | \le 2M.
    \]
    Utilizing Hoeffding's inequality,
\[
\begin{split}
    P_{IV} &\le P\left(\exists \theta\in\Theta_{i}, \frac{1}{n}\sum_{t=1}^{n} \left(g_{\theta}\left(\bX_{t}\right)-g_{s\left(\theta\right)}^{0,L}\left(\bX_{t}\right)\right)1\left(K_{s\left(\theta\right),t}=0\right)\ge \eta_{0}\right) \\
    & \le P\left(\exists \theta\in\Theta_{i}, \frac{1}{n}\sum_{t=1}^{n} \left(g_{\theta}^{0,U}\left(\bX_{t}\right)-g_{s\left(\theta\right)}^{0,L}\left(\bX_{t}\right)\right)1\left(K_{s\left(\theta\right),t}=0\right)\ge \eta_{0}\right) \\
    & \le N_{[]}\left(\mathcal{G}(\theta, r_{i}),\left\Vert \cdot\right\Vert_{L_2}, \delta_{0} \right) \exp\left(- \frac{n \eta_{0}^{2}}{4 M^2}\right) \\
    & \le \left( \frac{A r_i}{\delta_{0}} \right)^{m} \exp\left(- \frac{n \eta_{0}^{2}}{4 M^{2}}\right)
\end{split}
\]
Note that since 
\[
\frac{M \sqrt{2\ln (2)\frac{mi}{n}+ \frac{(i+1)(1-4\gamma)\xi}{n}}}{1/3 \gamma r_i} = \frac{M \sqrt{(1/2)\left(i+1\right)\left(1-4\gamma\right)\xi + (i+1)(1-4\gamma)\xi}}{1/3 \gamma 2^{i/2} \sqrt{\xi}} \underset{i \to \infty}{\to} 0,
\]
there exists a  constant \(C_1(\rho) > 0\) such that 
\[
\exp\left(- \frac{n \eta_{0}}{4 \ln^2(1/\delta)}\right) \le C_1(\rho) \exp\left(- \frac{\ln 2}{2} mi - \frac{(i+1)(1-4\gamma)\xi}{4}\right).
\]
Hence, 
\[
P_{IV} \le C_1(\rho) \left( \frac{A r_i}{\delta_{0}} \right)^{m} \exp\left(- \frac{\ln 2}{2} mi - \frac{(i+1)(1-4\gamma)\xi}{4}\right)
\]

\end{enumerate}

The rest of the proof is along the lines of the proof of Lemma 0 of Yang and Barron. Set
\[
\begin{split}
\delta_{j} &=Ar_{0}\exp\left(-\frac{\left(j+1\right)\left(1-4\gamma\right)\xi}{4m}\right) \\
\eta_{j} &=\delta_{j-1}\sqrt{8\ln\left(2\right)mi+\frac{2\left(2j+1\right)\left(1-4\gamma\right)\xi}{n}+\frac{\left(i+1\right)j\left(1-4\gamma\right)\xi}{n}}.
\end{split}
\]
 Note that for \(\rho\le A\), if \(\frac{\xi}{m}\ge\frac{4}{1-4\gamma}\ln\frac{2A}{\rho}\), then 
 \[
 \frac{\ln2}{2} mi\le\frac{\left(i+1\right)\left(1-4\gamma\right)\xi}{8}
 \]
 Since \(\eta_{j}\le A\sqrt{1-4\gamma}\sqrt{i+5}\frac{\xi}{n}\exp\left(-\frac{j\left(1-4\gamma\right)\xi}{8m}\right)\) (see Yang and Barron, 1998 for more details), we have 
 \[
 \begin{split}
 \sum_{j=1}^{\infty}\eta_{j}	&\le\sum_{j=1}^{\infty}A\sqrt{1-4\gamma}\sqrt{i+5}\frac{\xi}{n}\exp\left(-\frac{j\left(1-4\gamma\right)\xi}{8m}\right) \\
	&= A\sqrt{1-4\gamma}\sqrt{i+5}\frac{\xi}{n}\frac{\exp\left(-\frac{\left(1-4\gamma\right)\xi}{8m}\right)}{1-\exp\left(-\frac{\left(1-4\gamma\right)\xi}{8m}\right)} \\
	& \le A\sqrt{1-4\gamma}\sqrt{i+5}\frac{\xi}{n}\left(\frac{1}{1-\exp\left(-\frac{\left(1-4\gamma\right)\xi}{8m}\right)}\right)\exp\left(-\frac{\left(1-4\gamma\right)\xi}{8m}\right) \\
	& \le A\sqrt{1-4\gamma}\sqrt{i+5}\left(\frac{\sqrt{2}}{1-\sqrt{2}}\right)\frac{\xi}{n}\exp\left(-\frac{\left(1-4\gamma\right)\xi}{8m}\right) \\
	&\le2\sqrt{2}A\sqrt{1-4\gamma}\sqrt{i+5}\frac{\xi}{n}\exp\left(-\ln\frac{2A}{\rho}\right)
	=2\sqrt{2}\sqrt{1-4\gamma}\sqrt{i+5}\frac{\xi}{n}\rho
\end{split}
\]
Since \(\sqrt{i+5}\le\sqrt{5}2^{i}\), for \(\rho=\frac{\gamma}{6\sqrt{10}\sqrt{1-4\gamma}}\),
\[
\begin{split}
\sum_{j=1}^{\infty}\eta_{j}	&\le\frac{1}{3}\gamma2^{i}\frac{\xi}{n}=\frac{1}{3}\gamma r_{i}^{2} \\
\frac{\delta_{0}}{r_{0}} &\le A\exp\left(-\frac{\left(1-4\gamma\right)\xi}{4m}\right)=A\exp\left(-\ln\frac{2A}{\rho}\right)=\frac{\rho}{2}\le\rho.
\end{split}
\]
Hence, 
\[
\begin{split}
P_{I} &\le\exp\left(\frac{\ln2}{2}mi+\frac{\left(1-4\gamma\right)\xi}{4}\right)\exp\left(-\frac{\left(i+1\right)\left(1-4\gamma\right)\xi}{2}\right) \\
    &=\exp\left(\frac{\left(i+1\right)\left(1-4\gamma\right)}{8}-\frac{\left(i+1\right)\left(1-4\gamma\right)\xi}{4}\right) \\
	& \le\exp\left(-\frac{\left(i+1\right)\left(1-4\gamma\right)\xi}{8}\right) \\
P_{II},P_{III} &\le\sum_{j=1}^{\infty}\exp\left(\frac{\ln2}{2}mi-\frac{j\left(1-4\gamma\right)\xi}{4}\right)\exp\left(\frac{\ln2}{2}mi-\frac{\left(j+1\right)\left(1-4\gamma\right)\xi}{4}\right)\\
	&\cdot\exp\left(\ln\left(2\right)mi+\frac{\left(2j+1\right)\left(1-4\gamma\right)\xi}{4}-\frac{\left(i+1\right)j\left(1-4\gamma\right)\xi}{8}\right)\\
	&\le\sum_{j=1}^{\infty}\exp\left(-\frac{\left(i+1\right)j\left(1-4\gamma\right)\xi}{8}\right)
\end{split}
\]
\[
\begin{split}
    P_{IV} &\le C_1(\rho) \exp\left(\frac{\ln 2}{2} mi + \frac{(1-4\gamma)\xi}{4}\right) \exp\left(- \frac{\ln(2)}{2}mi - \frac{(i+1)(1-4\gamma)\xi}{4}\right) \\
    &\le C_1(\rho) \exp\left(\frac{\left(i+1\right)\left(1-4\gamma\right)}{8}-\frac{\left(i+1\right)\left(1-4\gamma\right)\xi}{4}\right) \\
	& \le C_1(\rho) \exp\left(-\frac{\left(i+1\right)\left(1-4\gamma\right)\xi}{8}\right) 
\end{split}
\]
Thus, 
\[
\begin{split}
P_{I}+P_{II}+P_{III} + P_{IV}&\le\left(1+\frac{2}{1-\exp\left(-\frac{\left(i+1\right)\left(1-4\gamma\right)\xi}{8}\right)} + C_1(M,\rho) \right)\exp\left(-\frac{\left(i+1\right)\left(1-4\gamma\right)\xi}{8}\right) \\
	&\le C(\rho) \exp\left(-\frac{\left(i+1\right)\left(1-4\gamma\right)\xi}{8}\right)
\end{split}
\]

\end{proof}